\documentclass[leqno]{article}
\usepackage{verbatim, amsmath, amscd, amssymb}
\usepackage[all]{xy}
\newcommand{\cA}{{\mathcal A}}
\newcommand{\cB}{{\mathcal B}}
\newcommand{\cC}{{\mathcal C}}

\newcommand{\cE}{{\mathcal E}}

\newcommand{\cL}{{\mathcal L}}

\newcommand{\cMod}{{\mathcal M}\!{\it od}}

\newcommand{\cO}{{\mathcal O}}
\newcommand{\cP}{{\mathcal P}}
\newcommand{\cQ}{{\mathcal Q}}

\newcommand{\rd}{\mathrm{d}}

\newcommand{\rD}{\mathrm{D}}

\newcommand{\rG}{\mathrm{G}}\newcommand{\rH}{\mathrm{H}}

\newcommand{\rN}{\mathrm{N}}

\newcommand{\rS}{\mathrm{S}}

\newcommand{\bbC}{\mathbb C}

\newcommand{\bbP}{\mathbb P}
\newcommand{\bbN}{\mathbb N}

\newcommand{\bbR}{\mathbb R}

\newcommand{\bbZ}{\mathbb Z}
\newcommand{\bfL}{\mathbf L}
\newcommand{\bfR}{\mathbf R}

\newcommand{\coh}{\mathrm{coh}}

\newcommand{\ch}{\mathrm{ch}\,}

\newcommand{\diag}{\mathrm{diag}}

\newcommand{\Dist}{\mathrm{Dist}}

\newcommand{\ev}{\mathrm{ev}}
\newcommand{\Ext}{\mathrm{Ext}}

\newcommand{\GL}{\mathrm{GL}}

\newcommand{\id}{\mathrm{id}}

\newcommand{\ind}{\mathrm{ind}}

\newcommand{\rad}{\mathrm{rad}}

\newcommand{\rk}{\mathrm{rk}\,}
\newcommand{\SL}{\mathrm{SL}}
\newcommand{\soc}{\mathrm{soc}}
\newcommand{\Sp}{\mathrm{Sp}}

\newcommand{\Mod}{\mathbf{Mod}}

\newcommand{\lbr}{\begin{bmatrix}}
\newcommand{\rbr}{\end{bmatrix}}
\newcommand{\for}{\bigcirc\kern-2.6ex \because}
\newcommand{\forb}{\bigcirc\kern-2.8ex \because}
\newcommand{\forbb}{\bigcirc\kern-3.0ex \because}
\newcommand{\forbbb}{\bigcirc\kern-3.1ex \because}

\newcommand\pf{\noindent {\bf Proof:  }}
\newtheorem{thm}{Theorem:}

\newtheorem{prop}{Proposition:}

\newtheorem{lem}{Lemma:}

\newtheorem{cor}{Corollary:}

\begin{document}
\large
\title{
\bf
Some observations 
on 
Karoubian 
complete strongly exceptional posets on the projective
\\
homogeneous varieties
\thanks
{The first author is supported in part by JSPS Grant in Aid
for Scientific Research, 
and the second author by The National Natural Science Foundation of China 10671142.
}
} 
\author{K\textsc{aneda} Masaharu\\558-8585 
Sugimoto\\
Osaka City University\\Graduate School of Science\\Department of Mathematics\\
kaneda@sci.osaka-cu.ac.jp\and
\and
Y\textsc{e} Jiachen
\\
Department of
Mathematics
\\
Tongji University
\\
1239 Siping Road
Shanghai 200092
\\
P. R. China
\\
jcye@mail.tongji.edu.cn}
\maketitle

\begin{abstract}
Let $\cP=G/P$ be a homogeneous projective variety with $G$ a reductive group and $P$ a parabolic subgroup.
In positive characteristic
we exhibit for $G$ of low rank  
a Karoubian complete strongly
exceptional poset of
locally free sheaves
appearing in the Frobenius direct image
of the structure sheaf of 
$G/P$.
These sheaves are all defined over
$\bbZ$, so by base change provide 
a Karoubian complete strongly
exceptional poset on
$\cP$ over
$\bbC$,
adding to the list of classical results by Beilinson and Kapranov
on the Grassmannians and the quadrics over $\bbC$.
\end{abstract}

On the complex projective space
$\bbP^n_\bbC$
Beilinson \cite{Bei}
found that
$\cE=\coprod_{i=0}^n\cO(-i)$
induces a triangulated
equivalence
from
the bounded derived category of
coherent sheaves on
$\bbP^n_\bbC$ to the bounded derived category of
right modules of finite type over the endomorphism ring of
$\cE$.
After the discovery of similar phenomena
by Kapranov
\cite{Kap}/
\cite{Kap88}
on the homogeneous
projective varieties
$\cP_\bbC=G_\bbC/P_\bbC$
with
$G_\bbC=\GL_n(\bbC)$ and a parabolic subgroup
$P_\bbC$, and on the complex quadrics,
Catanese \cite{Boh}
has proposed a conjecture on the existence of
coherent sheaves
$\cE_w$
on $\cP_\bbC$ for general complex
reductive group
$G_\bbC$
parametrized by
the coset representatives
of the Weyl group
of
$G_\bbC$ by the Weyl group
of
$P_\bbC$ such that
(i)
$\Mod_{\cP_\bbC}(\cE_w, \cE_w)\simeq\bbC$
$\forall w$,
(ii)
$\Ext^i_{\cP_\bbC}(\coprod_w\cE_w, \coprod_w\cE_w)=0$
$\forall i>0$,
(iii)
$\Mod_{\cP_\bbC}(\cE_x, \cE_y)\ne0$
iff
$x\leq y$
in the Chevalley-Bruhat order on $W$,
and 
(iv)
the 
$\cE_w$'s generate the bounded derived category of coherent sheaves
on
$\cP_\bbC$; precisely, we will consider Karoubian generation in this paper that the smallest triangulated subcategory
of
$\rD^b(\coh\cP)$ containing the 
$\cE_w$'s and closed under taking direct summands should be the whole of
$\rD^b(\coh\cP)$.
If the conjecture holds,
$\bfR\Mod_{\cP_\bbC}(\coprod_w\cE_w, ?)$
gives by Beilinson's lemma
\cite{Bei}/\cite{Ba}
a triangulated equivalence
from
the bounded derived category of
coherent sheaves on
$\cP_\bbC$ to that of
right modules of finite type over the endomorphism ring of
$\coprod_w\cE_w$,
and also supports Kontsevich's homological mirror conjecture
\cite{Boh}.

In this paper we propose a new way of prescribing
where to look for such
$\cE_w$'s.
We go over to
positive characteristic and
exhibit for $G$ of rank at most 2
the 
$\cE_w$'s 
as indecomposable direct summands
of the Frobenius direct image of the structure sheaf of
$\cP$;
in case $G$ is in type $\rG_2$
we do construct those
$\cE_w$'s based on such data 
but are unable at present to prove
that they indeed all appear in the Frobenius direct image.
To describe our method, we have to introduce some more notations.

Let
$\Bbbk$ be an algebraically closed field of
positive characteristic
$p$,
$G$ a reductive algebraic group over
$\Bbbk$, $P$ a parabolic subgroup
of $G$ and put
$\cP=G/P$.
We assume $p>h$ 
the Coxeter number of
$G$.
Let
$F:\cP\to\cP$ be the absolute Frobenius endomorphism of $\cP$.
If $G_1$ is the Frobenius kernel of
$G$,
$F$ factors through a natural morphism
$q:\cP\to G/G_1P$
to induce an isomorphism 
$\phi$
of schemes,
though not of
$\Bbbk$-schemes,
from
$G/G_1P$ to
$\cP$.
Now,
$q_*\cO_\cP\simeq\cL_{G/G_1P}(\hat
\nabla_P(\varepsilon))$
locally free $\cO_{G/G_1P}$-module
associated to $G_1P$-module
$\hat
\nabla_P(\varepsilon)$
induced from trivial 1-dimensional
$P$-module
$\varepsilon$
\cite{Haa}.
As
$F_*\cO_\cP\simeq\phi_*q_*\cO_\cP$,
the structure of
$F_*\cO_\cP$ is controlled by that of
$\hat
\nabla_P(\varepsilon)$.
Let $T$ be a maximal torus of $P$.
We propose a formula to describe the $G_1T$-socle series
on
$\hat
\nabla_P(\varepsilon)$ in terms of
Kazhdan-Lusztig polynomials
and examine
each socle layer, which is equipped with a structure of
$G_1P$-module.
If $W$ is the Weyl group of
$G$ with length function $\ell$ 
and $W_P$ the Weyl group of $P$,
$W^P=\{w\in W\mid
\ell(wx)=\ell(w)+\ell(x)\ \forall x\in W_P\}$
gives the set of
coset representatives of
$W/W_P$.
We find that the multiplicity space of
$G_1$-simple module parametrized by
$w\in W^P$ appearing in the
$(\ell(w)+1)$-st socle layer 
of $\hat
\nabla_P(\varepsilon)$
induce by sheafification the desired sheaf
$\cE_w$, inverting the order on
$W^P$.
Thus our investigation is twofold;
one is to study the structure of
induced
$G_1P$-modules,
and the other is to study the sheafification of those $G_1P$-modules
arising from the $G_1T$-socle series of
the induced module.

We note that our work is also related to the tilting property of
$F_*\cO_\cP$.
We observed in \cite{HKR} that
if 
$F_*\cO_\cP$ is tilting, 
the triangulated localization theorem holds for the endomorphism ring
of
$\cO_\cP$
over
its Frobenius twist
$\cO_\cP^{(1)}$,
which is a version of
Bezrukavnikov-Mirkovic-Rumynin
localization theorem for the $\Bbbk$-algebra  of crystaline differential operators on $\cP$
\cite{BMR}.
As
$\cO_\cP$ is locally free of rank
$p^{\dim\cP}$ over $\cO_\cP^{(1)}$,
the category of
coherent
$\cO_{\cP}$-modules is equivalent to
the category of
coherent modules over the sheaf of
small differntial operator ring
$\cMod_{\cO_\cP^{(1)}}(\cO_\cP,\cO_\cP)$.

Since the first author presented a talk on a part of the present work at Tongji
University in 2006, a number of related works have appeared.
In particular, we have verified in
\cite{K08}/\cite{KNS} that Kapranov's sheaves on the Grassmannian provide the desired sheaves
in positive characteristic,
while Langer \cite{La} has proved that
$F_*\cO_\cQ$ is tilting on the quadrics
$\cQ$,
see also Samokhin \cite{S07},
\cite{S1}-\cite{S3}.
One can parametrize certain direct summands of $F_*\cO_\cQ$
by $W^P$ to verify Catanese's conjecture on the quadrics.
On the projective spaces \cite{K09} has showed that 
$\hat\nabla_P(\varepsilon)$
is uniserial as
$G_1T$-module, and that the multiplicity
spaces of the $G_1$-simple
isotypic components in the
$G_1T$-socle layers in
$\hat\nabla_P(\varepsilon)$
provide the desired $\cE_w$'s
as proposed in this work.

Our main results are stated in \S1.
In \S\S2 and 3 we show that our $\cE_w$
possess the conjectured properties.
As the verifications are done by brute force, we omit tedious mechanical computations.
In \S4 we discuss parametrization of
Kapranov's sheaves based on our observations.

We are grateful to Henning Andersen, 
Hashimoto Yoshitake,
Tezuka Michishige and Yagita Nobuaki for helpful dicussions.
We learned of B\"{o}hning's preprint of
\cite{Boh} and also of
\cite{S07},
\cite{La} from Hashimoto.
Thanks are also due to Adrian Langer 
for explaining his work.
A part of the work was done
during the first author's visit to the second in Shanghai in the fall of 2006.
He thanks Tongji University for the hospitality and the financial support during the visit.
Some of the communications were made while the second author was visiting Abdus Salam International Centre for Theoretical Physics in the summer of 2007, to which
he thanks for 
the hospitality and the financial support.

\bigskip
\setcounter{equation}{0}
\begin{center}
$1^\circ$
{\bf
Parabolic Humphreys-Verma modules
}
\end{center}

In this section after fixing 
the notations to be used throughout the manuscript, we begin our study of induced $G_1P$-modules by relating them to better-examined
induced $G_1B$-modules, $B$ a Borel subgroup
of $G$.
We will determine the $G_1T$-socle series on parabolic Humphreys-Verma modules for $G$ of rank at most $2$.
Based on their structural data we will define
our sheaves 
$\cE_w$ on $G/P$,
which verify Catanese's conjecture.

\setcounter{equation}{0}
\noindent
(1.1)
We will assume $G$ is a simply connected 
simple algebraic group over
an algebraically closed field
$\Bbbk$ of characteristic
$p>0$.
Let $P$ be a parabolic subgroup of 
$G$,
$B$ a Borel subgroup of $P$, and
$T$ 
a maximal torus of $B$.
Let 
$\Lambda$ be the character group of
$B$,
$R\subseteq \Lambda$
the root system of
$G$ relative to $T$ with
positive system $R^+$ such that the roots of $B$ are $-R^+$.
If $\alpha\in R$, we denote its coroot by
$\alpha^\vee$.
Let $R^s$ be the set of simple roots of
$R^+$ 
and $\Lambda^+$ the corresponding
set of dominant weights.
If $\alpha\in R^s$, let
$\omega_\alpha\in\Lambda$ 
such that
$\langle\omega_\alpha,
\beta^\vee\rangle=\delta_{\alpha\beta}$
$\forall \beta\in
R^s$.
Let $W$
be the Weyl group of
$G$ with distinguished generators
$s_\alpha$,
$\alpha\in
R^s$.
For $w\in W$ and $\lambda\in\Lambda$ set
$w\bullet\lambda=w(\lambda+\rho)-\rho$
with
$\rho=
\frac{1}{2}\sum_{\alpha\in
R^+}\alpha$.
We will denote the Chevalley-Bruhat order (resp. the length function)
on $W$ relative to
$\{s_\alpha\mid
\alpha\in R^s\}$
by
$\geq$
(resp. $\ell$).
If 
$W^{P}=\{w\in W\mid \ell(wx)=\ell(w)+\ell(x) \ \forall x\in
W_P\}$, 
then
$W=\sqcup_{w\in W^{P}}wW_{P}$.
Let
$w_0$
(resp. $w_P$)
be the longest element of
$W$
(resp. $W_P$).
If we write $P=P_I$ for
$I\subset R^s$, let $R_I=R\cap\sum_{\alpha\in I}\bbZ\alpha$ the root system of the standard Levi subgroup of
$P$.

For an algebraic group $H$ over $\Bbbk$ let
$H_1$ be the Frobenius kernel of
$H$ and
$\Dist(H)$ the algebra of distributions of
$H$ over $\Bbbk$.
Let
$\hat\nabla_P=\ind_P^{G_1P}$
be the induction functor from the category of
$P$-modules to the category of $G_1P$-modules,
which we call the Humphreys-Verma induction.
Let
$\Lambda_{P}$ be the character group of
$P$. 
If $\nu\in\Lambda_P$
and if
$\rho_P=\frac{1}{2}\sum_{\alpha\in
R^+\setminus
R_I}\alpha$,
there is an isomorphism of
$G_1P$-modules
$
\hat\nabla_P(\nu)
\simeq
\Dist(G_1)\otimes_{\Dist(P_1)}(\nu-2(p-1)\rho_P)$.
Likewise we let
$\nabla_P=\ind_P^G$ be the induction from the category of $P$-modules to
the category of $G$-modules.
In case $P$ is $B$, we will often suppress
$B$ from the subscripts.
We also let
$\nabla^P=\ind_B^P$, abbreviated as
$\nabla^\alpha$ in case
$P=P_{\{\alpha\}}$ for a simple root $\alpha$.
If $\lambda\in\Lambda$,
we put
$\Delta^P(\lambda)=\bfR^{\dim P/B}\ind_B^P(w_P\bullet\lambda)$,
abbreviated as
$\Delta^\alpha(\lambda)$ in case
$P=P_{\{\alpha\}}$.
For an $H$-module $M$
we denote by
$M^*$ the
$\Bbbk$-linear dual of
$M$ and by $M^{[1]}$ the Frobenius twist of
$M$.
If $H_1$ acts trivially on $M$, one can untwist the Frobenius action and obtain an
$H$-module
$M^{[-1]}$ such that
$(M^{[-1]})^{[1]}\simeq M$.
By $\soc_{H}M$ we will denote the socle of
$M$ as $H$-module.

If $M$ is a $T$-module and if
$\nu\in\Lambda$,
$M_\nu$ will denote the $\nu$-weight space of
$M$.
Let
$\ch M=\sum_{\lambda\in\Lambda}(\dim M_\lambda)
e^\lambda$ be the formal character of
$M$.
The simple $G$-
(resp. $G_1T$-) modules
are parametrized by their highest weights in
$\Lambda^+$
(resp. $\Lambda$),
denoted
$L(\lambda)$,
$\lambda\in\Lambda^+$
(resp. $\hat L(\nu)$,
$\nu\in\Lambda$).
In order to avoid confusion,
we will let
$\varepsilon$ denote
trivial 1-dimensional $G$-module.
For other unexplained notations we refer to
\cite{J} except that
for a category
$\cC$ we will denote the set of morphisms of $\cC$ from object $A$ to $B$
by
$\cC(A,B)$.

\setcounter{equation}{0}
\noindent
(1.2)
If $I\subseteq R^s$ associated to
$P$,
$\Lambda_{P}=\{\lambda\in\Lambda\mid
\langle\lambda,\alpha^\vee\rangle=0
\ \forall\alpha\in I\}$.
Let us begin with relating
the formal character of
$\hat\nabla_P(\nu)$,
$\nu\in\Lambda_P$,
to that of $\hat\nabla(\lambda)$,
$\lambda\in\Lambda$.
As $G_1T$-module
$\hat\nabla_P(\nu)$ is isomorphic to
$G_1T$-module
$\ind_{P_1T}^{G_1T}(\nu)$
induced from $P_1T$-module
$\nu$.
Let
$\widehat{\bbZ[\Lambda]}$
be a completion of
$\bbZ[\Lambda]\subset
\Pi_{\lambda\in\Lambda}\bbZ e^\lambda$
as in \cite[3.3]{F}
consisting of
those
$(a_\lambda
e^\lambda)_\lambda
\in
\Pi_{\lambda\in\Lambda}\bbZ e^\lambda$
such that
there is
$\mu\in\Lambda$ for which
whenever
$a_\lambda\ne0$,
$\lambda\leq\mu$.

\begin{prop}
Let $L$ be the standard Levi subgroup of
$P$ 
and $U_L$ the unipotent radical
of the 
Borel subgroup
$B\cap L$ of
$L$.
$\forall\nu\in\Lambda_P$, 
one has in
$\widehat{\bbZ[\Lambda]}$
\[
\ch\hat\nabla_P(\nu)=
e^\nu\prod_{\alpha\in R^+\setminus R_I}
\frac{1-e^{-p\alpha}}{1-e^{-\alpha}}
=\sum_{ \substack{ w\in W_P
\\
\gamma\in\bbZ
R_I}}
(-1)^{\ell(w)}
\dim(\Dist(U_L)_\gamma)
\ch\hat\nabla(w\bullet\nu+p\gamma).
\]

\end{prop}

\pf 
The first equality follows from the decomposition
$G_1\simeq
P_1\times\prod_{\alpha\in R^+\setminus R_I}U_{\alpha,1}$
with
$U_{\alpha}$ root subgroup
associated to
$\alpha$.
If
$\Lambda_L$ is the character group of
$B\cap L$,
one can write in $\widehat{\bbZ[\Lambda_L]}$ 
\begin{align*}
\ch\varepsilon
&=
\sum_{\lambda\in\Lambda_L}
(\varepsilon : \ind_{(B\cap L)_1T}^{L_1T}(\lambda))
\ch\ind_{(B\cap L)_1T}^{L_1T}(\lambda)
\quad\text{for some
$(\varepsilon : \ind_{(B\cap L)_1T}^{L_1T}(\lambda))\in\bbZ$}
\\
&=
\sum_{\lambda\in\Lambda_L}
(\varepsilon : \ind_{B_1T}^{P_1T}(\lambda))
\ch\ind_{B_1T}^{P_1T}(\lambda).
\end{align*} 
Then in $\widehat{\bbZ[\Lambda]}$
\begin{align*}
\ch\hat\nabla_P(\nu)
&=
\ch\hat\nabla_P(\nu\otimes\varepsilon)
=
\sum_{\lambda\in\Lambda_L}
(\varepsilon : \ind_{B_1T}^{P_1T}(\lambda))
\ch\hat\nabla_P(\nu\otimes
\ind_{B_1T}^{P_1T}(\lambda))
\quad\text{as $\hat\nabla_P$ is exact}
\\
&=
\sum_{\lambda\in\Lambda_L}
(\varepsilon : \ind_{B_1T}^{P_1T}(\lambda))
\ch\hat\nabla_P(
\ind_{B_1T}^{P_1T}(\nu+\lambda))
\quad\text{by the tensor identity as
$\nu\in\Lambda_P$}
\\
&=
\sum_{\lambda\in\Lambda_L}
(\varepsilon : \ind_{B_1T}^{P_1T}(\lambda))
\ch\hat\nabla(\nu+\lambda)
\quad\text{by the transitivity of inductions}.
\end{align*}
On the other hand, 
Weyl's character formula for $L$ asserts
\begin{align*}
e^0
&=
\frac{\sum_{w\in W_P}(-1)^{\ell(w)}e^{w\bullet_L0}}
{\prod_{\alpha\in R_I^+}(1-e^{-\alpha})}
\quad\text{with
$w\bullet_L0=w\rho_L-\rho_L$,
$\rho_L=\frac{1}{2}\sum_{\alpha\in R_I^+}\alpha$}
\\
&=
\sum_{w\in W_P}(-1)^{\ell(w)}
e^{w\bullet_L0}
\prod_{\alpha\in R_I^+}\frac{1-e^{-p\alpha}}{1-e^{-\alpha}}
\prod_{\alpha\in R_I^+}\frac{1}{1-e^{-p\alpha}}
\\
&=
\prod_{\alpha\in R_I^+}\frac{1-e^{-p\alpha}}{1-e^{-\alpha}}
\sum_{\substack{w\in W_P
\\
\gamma\in\bbN
R_I^+}}(-1)^{\ell(w)}
\dim(\Dist(U_L)_{-\gamma})
e^{w\bullet_L0-p\gamma}
\\
&=
\sum_{\substack{w\in W_P
\\
\gamma\in\bbN
R_I^+}}(-1)^{\ell(w)}
\dim(\Dist(U_L)_{-\gamma})
\ch\,\ind_{(B\cap L)_1T}^{L_1T}(w\bullet_L0-p\gamma)
\\
&=
\sum_{\substack{w\in W_P
\\
\gamma\in\bbZ
R_I}}(-1)^{\ell(w)}
\dim(\Dist(U_L)_{\gamma})
\ch\,\ind_{B_1T}^{P_1T}(w\bullet_L0+p\gamma)
\\
&=
\sum_{\substack{w\in W_P
\\
\gamma\in\bbZ
R_I}}(-1)^{\ell(w)}
\dim(\Dist(U_L)_{\gamma})
\ch\,\ind_{B_1T}^{P_1T}(w\bullet0+p\gamma).
\end{align*}

It follows that
\begin{align*}
\ch\hat\nabla_P(\nu)
&=
\sum_{\substack{w\in W_P
\\
\gamma\in\bbZ
R_I}}(-1)^{\ell(w)}
\dim(\Dist(U_L)_{\gamma})
\ch\hat\nabla(\nu+w\bullet0+p\gamma)
\\
&=
\sum_{\substack{w\in W_P
\\
\gamma\in\bbZ
R_I}}(-1)^{\ell(w)}
\dim(\Dist(U_L)_{\gamma})
\ch\hat\nabla(w\bullet\nu+p\gamma)
\quad\text{as
$\langle\nu,\alpha^\vee\rangle=0$
$\forall\alpha\in
R_I$}.
\end{align*}

\setcounter{equation}{0}
\noindent
(1.3)
\noindent
{\bf Lemma:}
{\it
Let $\nu\in\Lambda_P$.

(i)
There are isomorphisms of
$G_1B$-modules
$
\soc_{G_1P}\hat\nabla_P(\nu)
\simeq
\hat L(\nu)
\simeq
\soc_{G_1B}\hat\nabla_P(\nu)$,
and hence one may regard
$\hat\nabla_P(\nu)$ as 
a $G_1B$-submodule of
$\hat\nabla(\nu)$.

(ii)
There is an isomorphism of
$G_1P$-modules
$
\hat\nabla_P(\nu)^*\simeq
\hat\nabla_P(2(p-1)\rho_P-\nu)$.

}

\pf
(i)
One has
\[
G_1B
\Mod(\ind_P^{G_1P}(\nu),
\ind_B^{G_1B}(\nu))
\simeq
B\Mod(\ind_P^{G_1P}(\nu),
\nu)
=
\Bbbk\,\ev_\nu.
\]
On the other hand,
$\forall\lambda\in\Lambda$,
\begin{align*}
G_1B
&\Mod(\hat L(\lambda),
\ind_P^{G_1P}(\nu))
\leq
G_1T\Mod(\hat L(\lambda),
\ind_{P_1T}^{G_1T}(\nu))
\simeq
P_1T\Mod(\hat L(\lambda),
\nu)
\\
&\leq
B_1T\Mod(\hat L(\lambda),
\nu)
=
G_1T\Mod(\hat L(\lambda),
\ind_{B_1T}^{G_1T}(\nu))
\simeq\Bbbk\delta_{\lambda\nu}.
\end{align*}
It follows that
$\soc_{G_1B}(\ind_P^{G_1P}(\nu))=\hat
L(\nu)$,
and hence also
$\soc_{G_1P}(\ind_P^{G_1P}(\nu))=\hat
L(\nu)$.
Thus
$\ind_P^{G_1P}(\nu)
\leq
\ind_B^{G_1B}(\nu)$
as $G_1B$-modules
via the commutative diagram
\[
\xymatrix{
\ind_P^{G_1P}(\nu)
\ar@{-->}[rr]
\ar[drr]_{\ev_\nu}
&&
\ind_B^{G_1B}(\nu)
\ar[d]^{\ev_\nu}
\\
&&
\nu.
}
\]

(ii) follows from
\cite[I.8.20
and II.3.4]{J}.

\setcounter{equation}{0}
\noindent
(1.4)
Assume now that
$p>h$ the Coxeter number of
$G$.
We say $\lambda\in\Lambda$ is $p$-regular iff
$\langle\lambda+\rho,\alpha^\vee\rangle\not\in p\bbZ$
$\forall\alpha\in R$.
Lusztig's conjecture on the irreducible characters for $G$, equivalently for $G_1T$,
is now a theorem for indefinitely large
$p$ due to
Andersen, Jantzen and Soergel
\cite{AJS}
and more recently to
Fiebig
\cite{F1}.
In turn, Lusztig's conjecture allows us to
determine the 
$G_1T$-socle series of
$\hat\nabla(\lambda)$ for $p$-regular
$\lambda$
\cite{AK89}.
To describe it, 
let
$W_p=W\ltimes p\bbZ R$ acting on $\Lambda$ with
$p\bbZ R$ by translations.
If $s_0$ is the reflexion with respect to the hyperplane
$\{v\in\Lambda\otimes_\bbZ\bbR\mid
\langle v+\rho, \alpha_0^\vee\rangle=-p\}$
with
$\alpha_0^\vee$ the highest coroot,
$(W_p, (s_{\alpha}, s_0\mid
\alpha\in R^s))$ forms a Coxeter system.
We consider also the translations by
$p\Lambda$.
Then
$\Lambda_1=\{\lambda\in\Lambda\mid
\langle\lambda,\alpha^\vee\rangle\in[0,p[\ \forall \alpha\in R^s\}$
is a fundamental domain for the action of
$p\Lambda$.
Put
$\cA=
(\Lambda\otimes_\bbZ\bbR)\setminus
\cup_{\alpha\in R, n\in\bbZ}\{v\in\Lambda\otimes_\bbZ\bbR\mid
\langle v+\rho,\alpha^\vee\rangle=pn\}$,
an element of which we call an alcove.
We will denote an alcove containing
$0$ by
$A^+$.
As the structure of $G_1T$-socle series is uniform for $p$-regular weights in an alcove,
let
$\hat\nabla(A)$
(resp. $\hat L(A)$), $A\in \cA$,
denote $\hat\nabla(\lambda)$
(resp. $\hat L(\lambda)$)
for $\lambda\in A$.
Let
$0<\soc^1\hat\nabla(A)=\soc\hat\nabla(A)\leq\soc^2\hat\nabla(A)\leq\dots$
be the $G_1T$-socle series of
$\hat\nabla(A)$
and let
$\soc_i\hat\nabla(A)=
\soc^i\hat\nabla(A)/\soc_{i-1}\hat\nabla(A)$
be the $i$-th socle layer of
$\hat\nabla(A)$.
The Lusztig conjecture implies that the
Loewy length, i.e., the length of the socle series of
$\hat\nabla(A)$ is
$\ell(w_0)+1$,
and that $\forall C\in\cA$,
\begin{equation}
Q^{C,A}
=
\sum_{i\in\bbN}q^{\frac{\rd(C,A)+1-i}{2}}
[\soc_{i}\hat\nabla(A) : \hat L(
C)],
\end{equation}
where
$Q^{C,A}$ is a Kazhdan-Lusztig polynomial in indeterminate
$q$
\cite[1.8]{L80}, 
$\rd(C,A)$
is the distance from alcove
$C$ to alcove
$A$,
and
$[\soc_{i}\hat\nabla(A) : \hat L(
C)]$
is the multiplicity
of
$\hat L(
C)$ as a $G_1T$-composition factor of
$\soc_i\hat\nabla(A)$.
For $G$ of rank $\leq2$ the formula 
(1) is known to hold for
$p\geq h$.

For
$A\in\cA$
let
$0_A$ be the element of
$W_p\bullet0$ in
$A$.
Based on (1), (1.2), and 
noting that
$
\soc_i\hat\nabla_P(A)
\leq
\soc_i\hat\nabla(A)$
$\forall i\in\bbN$,
we speculate
the
$G_1T$-socle series of
$\hat\nabla_P(A)$ with
$0_A\in\Lambda_P$
to be given by
\begin{multline}
\sum_{\substack{ w\in W_P
\\
\gamma\in\sum_{\alpha\in I}\bbZ\alpha}}
(-1)^{\ell(w)}
\dim(\Dist(U_L)_\gamma)
Q^{C, w\bullet A+{p\gamma}}
=
\sum_{i\in\bbN}
q^{\frac{d(C, A)+1-i}{2}}
[\soc_i\hat\nabla_P(A) : \hat L(C)].
\end{multline}
This is (1) in case $P=B$,
and specializes to
(1.2) under
$q\rightsquigarrow1$.
Put $w^P=w_0w_{P}$.

\noindent
{\bf Proposition:}
{\it
Assume $G$ is of rank $\leq2$.
For
$p\geq h$ 
each
$\hat\nabla_P(A)$ with $0_A\in\Lambda_P$ has
the Loewy length 
$\ell(w^P)+1$
and the $G_1T$-socle series is given by the formula (2).

}

\pf
We may assume $P>B$.
Unless $G$ is in type $\rG_2$,
$\hat\nabla(A)$ is multiplicity-free 
rendering the verification mechanical.

Assume now that $G$ is in type $G_2$.
Let
$\alpha_1$ and $\alpha_2$
be the simple roots with
$\alpha_1$ short, and put
$s_i=s_{\alpha_i}$ and $\omega_i=\omega_{\alpha_i}$.
We may translate $0_A$ into $\Lambda_1$ via $p\Lambda_P$. 
Let $P=P_{\alpha_1}$.
We may assume $A=A^+$ or $s_0s_1s_2s_1s_0\bullet A^+$.
Consider the case $A=A^+$.
We are to test for each $C\in\cA$
\begin{multline}
\sum_{
n\in\bbN}
Q^{C,A^+-pn\alpha_1}
-
\sum_{
n\in\bbN}
Q^{C,s_1\bullet A^+-pn\alpha_1}
=
\sum_i
q^{\frac{d(C,A^+)+1-i}{2}}
[\soc_i\hat\nabla_{P}(A^+) : \hat L(C)].
\end{multline}
If 
$Q^{C,A^+}$ is a monomial, 
we can read off the socle level of $\hat L(C)$ in $\hat \nabla_{P}(A^+)$ from (1.2), and (3) holds.
Thus those $C$ left to be examined are
$s_1s_2s_1s_0\bullet A^+-p\omega_2$,
$s_1s_2s_1s_0\bullet
A^+-2p\omega_1$,
$s_1s_0\bullet
A^+-p\omega_2$,
$A^++p(2\omega_1-2\omega_2)$,
$A^+-p\omega_2$,
$A^+-2p\omega_1$,
$s_1s_2s_1s_0\bullet
A^++p(\omega_1-2\omega_2)$,
$s_1s_2s_1s_0\bullet
A^+-p\rho$,
$A^++p(\omega_1-2\omega_2)$,
$A^+-p\rho$,
$A^+-3p\omega_1$,
$A^+-2p\omega_2$.
For (3) to hold, we must have all those
$\hat L(C)$ belonging to 
$\soc_{3}\hat\nabla_{P}(A^+)$.
As
$\hat L(s_1s_2s_1s_0\bullet
A^++p(2\omega_1-2\omega_2))=L(s_1s_2s_1s_0A^+)\otimes
(2\omega_1-2\omega_2)^{[1]}
\leq
\soc_{3}\hat\nabla_{P}(A^+)$
and as
$\Delta^{P}(2\omega_1-2\omega_2)$
is $P$-irreducible,
$G_1\Mod(L(s_1s_2s_1s_0\bullet A^+), \soc_3\hat\nabla_{P}(A^+))^{[-1]}\geq
\Delta^{P}
(2\omega_1-2\omega_2)$.
It follows that both
$\hat L(s_1s_2s_1s_0\bullet
A^+-p\omega_2)$
and
$\hat L(s_1s_2s_1s_0\bullet A^+-2p\omega_1)$
must also belong to
$\soc_3\hat\nabla_{P}(A^+)$.
Likewise
$\hat L(A^++p(\omega_1-2\omega_2))$,
$\hat L(A^+-p\rho)$,
$\hat L(A^+-3p\omega_1)$.
We will then be left with
the following
$C$:
$s_1s_0\bullet
A^+-p\omega_2$,
$A^++p(2\omega_1-2\omega_2)$,
$A^+-p\omega_2$,
$A^+-2p\omega_1$,
$s_1s_2s_1s_0\bullet
A^++p(\omega_1-2\omega_2)$,
$s_1s_2s_1s_0\bullet
A^+-p\rho$,
$A^+-2p\omega_2$.
For those $C$ if $\hat L(C)$ does not belong to
the third socle layer, it must lie in the 5th by (1).
On the other hand, dualizing an exact sequence of
$G_1P$-modules
\[
0\to
\soc_5\hat\nabla_{P}(A^+)
\to
\hat\nabla_{P}(A^+)/
\soc^4\hat\nabla_{P}(A^+)
\to
\hat L(s_0s_1s_2s_1s_0\bullet
A^+-2p\omega_2)
\to
0,
\]
one obtains another exact sequence
\begin{multline*}
0\to
\hat L(s_0s_1s_2s_1s_0\bullet
A^+-2p\omega_2)^*
\to
(\hat\nabla_{P}(A^+)/
\soc^4\hat\nabla_{P}(A^+)
)^*
\to
(\soc_5\hat\nabla_{P}(A^+))^*
\to
0
\end{multline*}
with
$
(\hat\nabla_{P}(A^+)
/
\soc^4\hat\nabla_{P}(A^+)
)^*
\leq
\hat\nabla_{P}(A^+)^*
\simeq
\hat\nabla_{P}(3(p-1)\omega_2)
\leq
\hat\nabla(3(p-1)\omega_2)
$
by (1.3), and hence
if $\hat L(C)\leq
\soc_5\hat\nabla_{P}(A^+)$,
then 
$
\hat L(C)^*\leq
\soc_2\hat\nabla(3(p-1)\omega_2)
=
\soc_2\hat\nabla(s_0s_1s_2s_1s_0\bullet
A^++2p\omega_2),
$
which contradicts (1).
Thus (3) holds.

Likewise the other cases.

\setcounter{equation}{0}
\noindent
(1.5)
Assume $p>h$.
For each $w\in W$
there is a unique element in
$(w\bullet0+p\Lambda)\cap\Lambda_1$,
which we will denote by
$\varepsilon_w$.
Thus the principal $G_1$-block consists of
$L(\varepsilon_w)$,
$w\in W$.
Put for simplicity
$L(w)=L(\varepsilon_w)$.
We know from \cite{Y} that
all $L(w)$, $w\in W$, appear
as
$G_1$-composition factors of
$\hat\nabla(\varepsilon)$.

\begin{cor}
Assume $\rk G\leq2$.
\begin{enumerate}
\item
Each $i$-th $G_1T$-socle layer of
$\hat\nabla(\varepsilon)$
admits a decomposition as
$G_1P$-module
\[
\soc_i\hat\nabla(\varepsilon)=
\coprod_{w\in W^P}L(w)\otimes
G_1\Mod(L(w),\soc_i\hat\nabla(\varepsilon)).
\]

\item
Each $L(w)$, $w\in W^P$, appears as $G_1$-factor of
$\soc_{\ell(w)+1}\hat\nabla_P(\varepsilon)$,
i.e.,
\linebreak
$G_1\Mod(L(w),\soc_{\ell(w)+1}\hat\nabla(\varepsilon))\ne0$.

\end{enumerate}

\end{cor}

\setcounter{equation}{0}
\noindent
(1.6)
{\bf Remarks:}
(i)
In case
$\cP=\GL(E)/P$
with
$\Bbbk$-linear space $E$ of basis 
$e_1,e_2,\dots, e_{n+1}$
and
with
$P=N_{\GL(E)}(\Bbbk e_{n+1})$,
the assertions (i) and (ii) hold
\cite{K09}.

(ii)
Over $\bbC$ if
$P^+_\bbC$  is the parabolic subgroup
opposite to
$P_\bbC$,
the $\Dist(G_\bbC)$-composition factors
of Verma module
$\Dist(G_\bbC)\otimes_{\Dist(P^+_\bbC)}\varepsilon_\bbC$
are known to be of the form
$L(w^{-1}\bullet0)$,
$w\in W^P$
\cite{Hum}.

\setcounter{equation}{0}
\noindent
(1.7)
Assume $\rk G\leq2$ and $p>h$.
Put
$\soc_{i,w}^1=G_1\Mod(L(w),\soc_i\hat\nabla(\varepsilon))^{[-1]}$
$\forall i\in[1, \ell(w^P)+1], w\in
W^P$.
Sheafifying
$\hat\nabla_P(\varepsilon)$
one obtains from (1.6)
a filtration on
$F_*\cO_\cP\simeq\phi_*\cL_{G/G_1P}(\hat\nabla_P(\varepsilon))$
of subquotients
$L(w)\otimes_\Bbbk\cL_\cP(\soc_{i,w}^1)$,
$i\in[1, \ell(w^P)+1], w\in
W^P$.
For our second objective it is therefore
important to determine
the $P$-module structure on
each
$\soc_{\ell(w)+1,w}^1$,
$w\in W^P$.
In case $G$ has rank $2$, let
$\alpha_1$ and $\alpha_2$ be the simple roots with
$\alpha_1$ short.
Put
$\omega_i=\omega_{\alpha_i}$ and
$s_i=s_{\alpha_i}$,
$i=1,2$.
Let $P_{\alpha_i}$ be the standard parabolic subgroup of $G$ associated to
$\alpha_i$, i.e., such that
$\pm\alpha_i$ are roots of
$P_{\alpha_i}$.
Arguing as in \cite{AK00}/\cite{HKR}/\cite{KY} we find
\begin{prop}
We have the following identifications as
$P$-modules.
\begin{enumerate}
\item
If $G=\SL_2$, then
$\soc_{1,e}^1\simeq\varepsilon$
and
$\soc_{2,w_0}^1\simeq-\rho$.

\item
If $G=\SL_3$ and $P=P_{\alpha_1}$,
then
$
\soc_{1,e}^1\simeq\varepsilon$,
$
\soc_{2,s_2}^1=-\omega_2$
and
$
\soc_{3,w^P}^1=-2\omega_2$.

\item
If $G=\SL_3$ and
$P=B$,
\begin{alignat*}{3}
\soc_{1,e}^1
&\simeq
\varepsilon,
&
\soc_{2,s_1}^1
&\simeq
-\omega_1
&
\soc_{2,s_2}^1
&\simeq
-\omega_2,
\\
\soc_{3,s_1s_2}^1
&\simeq
(-\rho)\otimes_\Bbbk
\Delta^{\alpha_2}(\omega_2),
\quad&
\soc_{3,s_2s_1}^1
&\simeq
(-\rho)\otimes_\Bbbk
\Delta^{\alpha_1}(\omega_1),
\quad&
\soc_{4,w_0}^1
&\simeq
-\rho.
\end{alignat*}

\item
If
$G=\Sp_4$ and $P=P_{\alpha_2}$,
\[
\soc_{1,e}^1
\simeq
\varepsilon,
\quad
\soc_{2,s_1}^1
\simeq
-\omega_1,
\quad 
\soc_{3,s_2s_1}^1
\simeq
-2\omega_1,
\quad 
\soc_{4,w^P}^1
\simeq
-3\omega_1.
\]

\item
If
$G=\Sp_4$ and $P=P_{\alpha_1}$,
\[
\soc_{1,e}^1
\simeq
\varepsilon,
\quad
\soc_{2,s_2}^1
\simeq
-\omega_2,
\quad
\soc_{3,s_1s_2}^1
\simeq
\Delta^{\alpha_1}(\omega_1-2\omega_2),
\quad
\soc_{4,w^P}^1
\simeq
-2\omega_2.
\]

\item
If $G=\Sp_4$ and $P=B$,
\begin{alignat*}{2}
\soc_{1,e}^1
&\simeq
\varepsilon,
&
\soc_{2,s_1}^1
&\simeq
-\omega_1,
\\
\soc_{2,s_2}^1
&\simeq
-\omega_2,
&
\soc_{3,s_1s_2}^1
&\simeq
(-\omega_2)\otimes_\Bbbk
\ker(\Delta(\omega_1)\twoheadrightarrow\omega_1),
\\
\soc_{3,s_2s_1}^1
&\simeq
(-\rho)\otimes
\Delta^{\alpha_1}(\omega_1),
&
\soc_{4,s_1s_2s_1}^1
&\simeq
(-\rho)\otimes_\Bbbk(\Delta(\omega_2)/(-\omega_2)),
\\
\soc_{4,s_2s_1s_2}^1
&\simeq
(-\rho)\otimes_\Bbbk(\Delta(\omega_1)/(-\omega_1)),
&
\soc_{5,w_0}^1
&\simeq
-\rho.
\end{alignat*}

\item
If $G$ is of type $\rG_2$ and $P=P_{\alpha_2}$,
\begin{align*}
\soc_{1,e}^1
&\simeq
\varepsilon,
\qquad
\soc_{2,s_1}^1
\simeq
-\omega_1,
\qquad
\soc_{3,s_2s_1}^1
\simeq
-2\omega_1,
\\
\soc_{4,s_1s_2s_1}^1
&\simeq
(-3\omega_1)\otimes
(\Delta(\omega_1)/\Dist(P_{\alpha_2})
(\Delta(\omega_1)_{-\alpha_1})),
\\
\soc_{5,s_2s_1s_2s_1}^1
&\simeq
-3\omega_1,
\qquad
\soc_{6,w^P}^1
\simeq
-4\omega_1.
\end{align*}

\item
If $G$ is of type $\rG_2$ and $P=P_{\alpha_1}$,
\begin{align*}
\soc_{1,e}^1
&\simeq
\varepsilon,
\qquad
\soc_{2,s_2}^1
\simeq
-\omega_2,
\\
\soc_{3,s_1s_2}^1
&\simeq
(-\omega_2)\otimes\ker(\Delta(\omega_1)\twoheadrightarrow
\Delta^{\alpha_1}(\omega_1)),
\qquad
\soc_{4,s_2s_1s_2}^1
\simeq
\Delta^{\alpha_1}(\omega_1-2\omega_2),
\\
\soc_{5,s_1s_2s_1s_2}^1
&\simeq
(-2\omega_2)\otimes(\Delta(\omega_1)/\Delta^{\alpha_1}(\omega_1-\omega_2)),
\qquad
\soc_{6,w^P}^1
\simeq
-2\omega_2.
\end{align*}

\end{enumerate}

\end{prop}

\setcounter{equation}{0}
\noindent
(1.8)
If $G$ is in type $\rG_2$ and $P=B$, we are unfortunately not able to detemine the $B$-module structure on
$\soc_{i,w}^1$ at present.
To speculate, we make use of the following

\begin{lem}
Let $G$ be an arbitrary
simply connencted
reducive algebraic group over
$\Bbbk$.
Assume
$p\geq2(h-1)$ and 
that
Lusztig's conjecture on
irreducible
$G_1T$-modules hold.
Let
$\rad_{G,i}\nabla(p\rho)$,
$i\in\bbN$,
be the $i$-th layer of
the radical series on
$\nabla(p\rho)$
as $G$-module,
and
for
each
$\lambda\in\Lambda_1$
let
$L(\lambda)\otimes(\rad_{G,i,\lambda}^1\nabla(p\rho))^{[1]}$
(resp.
$L(\lambda)\otimes(\soc_{i}^1\hat\nabla(\varepsilon))^{[1]}$)
denote the
$L(\lambda)$-isotypic component
of
$\rad_{G,i}\nabla(p\rho)$
(resp. $\soc_{i}\hat\nabla(\varepsilon)$
as $G_1T$-module).
Then there is for each $i$
a surjective homomorphism of
$B$-modules 
\[
(-\rho)\otimes_\Bbbk
\rad_{G,i,\lambda}^1\nabla(p\rho)
\twoheadrightarrow
\soc_{\ell(w_0)+1-i,\lambda}^1\hat\nabla(\varepsilon).\]

\end{lem}

\pf
By the hypothesis on $p$
one has
$\langle\nu^1+\rho,
\alpha^\vee\rangle\leq p$
for any weight $\nu$
of
$\nabla(p\rho)$
and any $\alpha\in R^+$,
where
$\nu=\nu^0+p\nu^1$
with
$\nu^0\in\Lambda_1$
and
$\nu^1\in\Lambda$.
Then by
\cite[8.1.2]{AK89}
the radical series as
$G$-module
and as $G_1T$-module on
$\nabla(p\rho)$
coincide:
$\forall i\in\bbN$,
\begin{equation}
\rad_{G,i}\nabla(p\rho)
=
\rad_{G_1T,i}\nabla(p\rho).
\end{equation}
On the other hand,
as observed in the proof of
\cite[8.2]{AK89},
the natural
homomorphism 
$
\nabla(p\rho)\to
\hat\nabla(p\rho)
$
of
$G_1B$-modules
is surjective,
which therefore induces
by
(1) 
a commutative diagram
of
$G_1T$-modules
\[
\xymatrix@R2ex{
\rad_{G,i}\nabla(p\rho)
\ar@{->>}[rr]
\ar@{-->>}
[rrdd]!<-10ex,2ex> 
&&
\rad_{G_1T,i}\hat\nabla(p\rho)
\ar@{=}[d]
\\
&&
\soc_{\ell(w_0)+1-i}\hat\nabla(p\rho)
\ar@{-}[d]^\sim
\\
&&
p\rho\otimes_\Bbbk
\soc_{\ell(w_0)+1-i}\hat\nabla(\varepsilon).
}
\]
As
$\rad_{G_1T,i}\hat\nabla(p\rho)$
is equipped with a structure of
$G_1B$-module,
the induced map 
\linebreak
$\rad_{G,i}\nabla(p\rho)\to
p\rho\otimes_\Bbbk
\soc_{\ell(w_0)+1-i}\hat\nabla(\varepsilon)$
above
is a surjective homomorphism of
$G_1B$-modules,
and hence the assertion. 

\setcounter{equation}{1}
\noindent
(1.9)
Assuming the Jantzen conjecture on the filtration
of
Weyl modules
\cite{A83}, the radical series of
$\Delta(p\rho)$ as $G$-module is available from
\cite{A87}.
Based also on the character data (1.4),
we are thus led to speculate on the $B$-module structure of $\soc_{\ell(w)+1,w}^1$,
$w\in W$,
in type
$\rG_2$
to be as follows:
\begin{align*}
\soc_{1,e}^1
&\simeq\varepsilon,
\\
\soc_{2,s_1}^1
&\simeq
-\omega_1,
\\
\soc_{2,s_2}^1
&\simeq
-\omega_2,
\\
\soc_{3,s_2s_1}^1
&\simeq
(-\rho)\otimes\Delta^{\alpha_1}(\omega_1),
\\
\soc_{3,s_1s_2}^1
&\simeq
(-\omega_2)\otimes\ker(\Delta(\omega_1)\twoheadrightarrow\omega_1),
\\
\soc_{4,s_2s_1s_2}^1
&\simeq
(-\rho)\otimes(\Delta(\omega_1)/\Dist(P_{\alpha_2})v_3),
\\
\soc_{4,s_1s_2s_1}^1
&\simeq
(-\rho)\otimes\{(\Delta(\omega_2)\oplus\Delta(\omega_1))/\Dist(P_{\alpha_1})(\Bbbk(v_2+v_1)+\Delta(\omega_2)_{-\alpha_2})\},
\\
\soc_{5,s_1s_2s_1s_2}^1
&\simeq
(-\rho)\otimes\{(\Delta(\omega_2)\oplus\Delta(\omega_1))/\Dist(P_{\alpha_2})(\Bbbk(v_4+v_3)+\Delta(\omega_2)_{-3\omega_1+\omega_2})\},
\\
\soc_{5,s_2s_1s_2s_1}^1
&\simeq
(-\rho)\otimes(\Delta(\omega_1)/\Delta^{\alpha_1}(\omega_1-\omega_2)),
\\
\soc_{6,s_2s_1s_2s_1s_2}^1
&\simeq
(-\rho)\otimes(\Delta(\omega_1)/(-\omega_1)),
\\
\soc_{6,s_1s_2s_1s_2s_1}^1
&\simeq
(-\rho)\otimes\{\varepsilon\oplus(\Delta(\omega_2)/(-\omega_2))\},
\\
\soc_{7,w_0}^1
&\simeq
-\rho,
\end{align*}
where
$v_i\in
\Delta(\omega_i)_{\alpha_1}\setminus0$,
$i\in\{1,2\}$,
and 
$v_{i+2}=\displaystyle
F_1^{(2)}v_i
\in
\Delta(\omega_i)_{-\alpha_1}$,
$i\in\{1,2\}$,
with
$F_1$ a root vector belonging to
$-\alpha_1$
in the Chevalley basis of the Lie algebra
of
$G$.

We remark, in particular, that
the Jantzen conjecture implies 
in type
$\rG_2$
that
$\soc^1_{6,s_1s_2s_1s_2s_1}$ 
should be decomposable, contrary to 
the cases for $\SL_2$, $\SL_3$, $\Sp_4$.

\setcounter{equation}{1}
\noindent
(1.10)
We now set
$\cE_w=\cL_\cP(\soc_{\ell(w)+1,w}^1)$,
$w\in W^P$,
or rather
we actually employ for
$(\cE_w)_{P/P}\otimes_{\cO_{\cP,P/P}}\Bbbk$
the $P$-module on the right hand side
of 
$\soc_{\ell(w)+1, w}^1$
given in (1.7) and (1.9), and 
forget about the characteristic restriction.
In type
$\rG_2$ when $P=B$
we take
$\cE_{s_1s_2s_1s_2s_1}=
\cL_{\cP}((-\rho)\otimes(\Delta(\omega_2)/(-\omega_2)))$
instead.

\begin{thm}
Let
$p>0$ be arbitrary unless otherwise specified.
\begin{enumerate}
\item
The
$\cE_w$, $w\in W^P$, Karoubian generate
$\rD^b(\coh\cP)$.

\item
Unless $G$ is in type $\rG_2$ with
$P=P_{\alpha_1}$ or $B$,
$\Mod_\cP(\cE_w, \cE_w)\simeq\Bbbk$
$\forall w\in W^P$,
$\Mod_\cP(\cE_x, \cE_y)\ne0$
iff
$x>y$
$\forall x,y\in W^P$,
and
$\Ext_\cP^i(\coprod_{w\in W^P}\cE_w, \coprod_{w\in W^P}\cE_w)\ne0$
$\forall i>0$.

\item
If $G$ is in type $\rG_2$ and
$P=P_{\alpha_1}$, assume $p\geq3$.
Then the same hold true for the
$\cE_w$, $w\in W^P$, as in (ii).

\item
If $G$ is in type $\rG_2$ and
$P=B$, assume $p\geq7$.
Then the same hold true for the
$\cE_w$, $w\in W^P=W$, as in (ii).

\end{enumerate}

\end{thm}

The arguments for the theorem
will be given in \S\S 2 and 3.

\setcounter{equation}{1}
\noindent
(1.11)
{\bf Remarks:}
(i)
In (iii) above the restriction on $p$ is sharp.
If $p=2$, we have,
see
\S3, $\forall i\in\bbN$,
\[
\Ext^i_\cP
(\cE_{s_2},
\cE_{s_1s_2})
\simeq
\begin{cases}
\Bbbk
&\text{if $i=0, 1$}
\\
0
&\text{else}.
\end{cases}
\]

(ii)
In (iv) if $p=2$,
we have $\forall i\in\bbN$,
\[
\Ext_\cP^i(\cE_{s_2}, \cE_{s_1s_2})
\simeq
\begin{cases}
\Bbbk
&
\text{if $i=0, 1$},
\\
0
&
\text{else}.
\end{cases}
\]
Also if $p=3$, $\Mod_\cP(\cE_{s_2s_1},
\cE_{s_1s_2s_1})\ne0$.

(iii)
If we employ the order reversing
involution
$w\mapsto w_0ww_P$ on
$W^P$, Theorem 1.10
verifies 
for groups of rank $\leq2$
Catanese's conjecture over
$\bbC$ by base change.
Our present parametrization 
appears more natural particularly
in the case of quadrics, see \S4.

\setcounter{equation}{0}
\begin{center}
$2^\circ$
{\bf
Karoubian completeness
}
\end{center}

Let 
$\{
\cE_w|
w\in W^P\}$
be the
coherent sheaves on $\cP=G/P$
defined in
(1.7) and (1.9).
In this section we will show
that
they are
Karoubian complete for
$\rD^b(\coh\cP)$, i.e.,
they Karoubian generate 
the bounded derived category
$\rD^b(\coh\cP)$
of coherent sheaves on
$\cP$.

Write
$P=P_I$, $I\subset R^s$,
and recall
$2\rho_P=
\sum_{\alpha\in
R^+\setminus
R_I}\alpha$.
As
$\cL_\cP(2\rho_P)$
is ample on $\cP$,
to show that the $\cE_w$,
$w\in
W^P$,
Karoubian generate $\rD^b(\coh\cP)$,
it is enough by a result
attributed to Kontsevich
by Positselskii
\cite{BMR}
to verify that all
$\cL_\cP(-2n\rho_P)$, $n\in\bbN^+$,
are Karoubian generated by the 
$\cE_w$'s; we actually show that
all
$\cL_\cP(\lambda)$,
$\lambda\in\Lambda_P$,
are Karoubian generated by the 
$\cE_w$'s.
This has been done in the case of 
$G\in\{
\SL_2, \SL_3\}$
with
$P=B$ in
\cite{HKR}
and
in the case of
$G=\Sp_4$ with $P=B$ in
\cite{KY};
precisely, in \cite{HKR} and in \cite{KY}
it was assumed that
$p\geq h$,
whose arguments carry over,
however, verbatim to
arbitrary characteristic
thanks to the new definition of
the
$\cE_w$'s.
Our argument
for $G$ in type $\rG_2$
and $P=B$ is essentially the same.

If $P>B$ in cases $G=\Sp_4$ or in type
$\rG_2$, we will also make use of a projection formula
\cite[3.3.2]{Boh}: put
$\cB=G/B$.
If
$\bar\pi : \cB\to\cP$ is the natural morphism,
\begin{equation}
\id_{\rD^b(\coh\cP)}
\simeq
(\bfR\bar\pi_*\cO_\cB)\otimes^{\bfL}_\cP\ ?
\simeq
(\bfR\bar\pi_*)\circ\bar\pi^* :
\rD^b(\coh\cP)
\to
\rD^b(\coh\cP).
\end{equation}
Recall from
\cite[I.5.19]{J}
that
$\forall m,n\in\bbZ$,
\begin{equation}
(\bfR\pi_*)
\cL_\cB(m\omega_1+n\omega_2)
\simeq
\begin{cases}
\cL_\cP(\nabla^{\alpha_1}(m\omega_1+n\omega_2))
\simeq
(\bfR\bar\pi_*)\cL_\cB(\nabla^{\alpha_1}(m\omega_1+n\omega_2))
\\
\hspace{5cm}
\text{if $P=P_{\alpha_1}$
and $m\geq0$}
\\
\cL_\cP(\nabla^{\alpha_2}(m\omega_1+n\omega_2))
\simeq
(\bfR\bar\pi_*)\cL_\cB(\nabla^{\alpha_2}(m\omega_1+n\omega_2))
\\
\hspace{5cm}
\text{if $P=P_{\alpha_2}$
and $n\geq0$},
\end{cases}
\end{equation}
and from
\cite[I.5.17]{J} that
$\forall M\in
P\Mod$,
\begin{equation}
\bar\pi^*\cL_\cP(M)
\simeq
\cL_\cB(M).
\end{equation}

In the following
we will often describe the setup in which a module $M$ admits a filtration whose subquotients are $M_1$,
$M_2$,\dots, $M_r$ from the bottom to the top by a diagram
\[
M=
\begin{tabular}{|c|}
\hline
$M_r$
\\
\hline
$\vdots$
\\
\hline
$M_2$
\\
\hline
$M_1$
\\
\hline
\end{tabular}.
\]
If $M_i$ is a direct sum of
$M_{i1},\dots,M_{is}$,
we will insert
\begin{tabular}{|c|}
\hline
$M_{i1}\mid
M_{i2}\mid\dots\mid
M_{is}$
\\
\hline
\end{tabular}
in place of
$M_i$.

Let us illustrate our argument
in case $G$ is in type
$\rG_2$ and $P=P_{\alpha_2}$.
Let
$\hat\cE=\langle\cE_w
\mid
w\in
W^P\rangle$ denote the triangulated subcategory of
$\rD^b(\coh\cP)$
Karoubian generated by
$\cE_w$,
$w\in
W^P$.
One has
$\Lambda_P=\bbZ\omega_1$.
Note first that
in
$\cE_{s_1s_2s_1}
=
\cL_{\cP}
((-3\omega_1)\otimes
(\Delta(\omega_1)/
\Dist(P)(v_3)
))$ one has
\[
\Dist(P)(v_3)
=
\begin{tabular}{|c|}
\hline
$\Delta^{\alpha_2}(-2\omega_1+\omega_2)$
\\ 
\hline
$-\omega_1$
\\ 
\hline
\end{tabular},
\]
and hence
\begin{align*}
(-3\omega_1)\otimes
(\Delta(\omega_1)/
\Dist(P)(v_3))
&\simeq
(-3\omega_1)\otimes\begin{tabular}{|c|}
\hline
$\omega_1$
\\
\hline
$\Delta^{\alpha_2}(-\omega_1+\omega_2)$
\\ 
\hline
$\varepsilon$
\\ 
\hline
\end{tabular}
=
\begin{tabular}{|c|}
\hline
$-2\omega_1$
\\
\hline
$\Delta^{\alpha_2}(-4\omega_1+\omega_2)$
\\ 
\hline
$-3\omega_1$
\\ 
\hline
\end{tabular}.
\end{align*}
It follows, by the presence of
$\cE_{s_2s_1s_2s_1}=\cL_\cP(-3\omega_1)$
and
$\cE_{s_2s_1}=\cL_\cP(-2\omega_1)$
in
$\hat\cE$,
that
$\cL_\cP(\nabla^{\alpha_2}(-4\omega_1+\omega_2))
\simeq
\cL_\cP(\Delta^{\alpha_2}(-4\omega_1+\omega_2))\in\hat\cE$.
Then, as
\begin{align*}
\hat\cE
&\ni\Delta(\omega_1)\otimes_\Bbbk
\cL_\cP(-2\omega_1)
\\
&\simeq
\cL_\cP(\Delta(\omega_1)\otimes_\Bbbk
(-2\omega_1))
\quad\text{by the tensor identity}
\\
&\simeq
\cL_\cP(\begin{tabular}{|c|}
\hline
$\omega_1$
\\
\hline
$\Delta^{\alpha_2}(-\omega_1+\omega_2)$
\\
\hline
$\varepsilon$
\\
\hline
$\Delta^{\alpha_2}(-2\omega_1+\omega_2)$
\\
\hline
$-\omega_1$
\\
\hline
\end{tabular})
\otimes
(-2\omega_1))
\quad
\text{using the 
$P$-module structure on
$\Delta(\omega_1)$}
\\
&\simeq
\cL_\cP(\begin{tabular}{|c|}
\hline
$-\omega_1$
\\
\hline
$\Delta^{\alpha_2}(-3\omega_1+\omega_2)$
\\
\hline
$-2\omega_1$
\\
\hline
$\Delta^{\alpha_2}(-4\omega_1+\omega_2)$
\\
\hline
$-3\omega_1$
\\
\hline
\end{tabular}),
\end{align*}
one obtains
$
\cL_\cP(\nabla^{\alpha_2}(-3\omega_1+\omega_2))
\simeq
\cL_\cP(\Delta^{\alpha_2}(-3\omega_1+\omega_2))
\in\hat\cE$.
In turn, as
\begin{align*}
\hat\cE
&\ni\Delta(\omega_1)
\otimes_\Bbbk\cL_\cP(-\omega_1)
\simeq
\cL_\cP(\begin{tabular}{|c|}
\hline
$\omega_1$
\\
\hline
$\Delta^{\alpha_2}(-\omega_1+\omega_2)$
\\
\hline
$\varepsilon$
\\
\hline
$\Delta^{\alpha_2}(-2\omega_1+\omega_2)$
\\
\hline
$-\omega_1$
\\
\hline
\end{tabular})
\otimes
(-\omega_1))
\simeq
\cL_\cP(\begin{tabular}{|c|}
\hline
$\varepsilon$
\\
\hline
$\Delta^{\alpha_2}(-2\omega_1+\omega_2)$
\\
\hline
$-\omega_1$
\\
\hline
$\Delta^{\alpha_2}(-3\omega_1+\omega_2)$
\\
\hline
$-2\omega_1$
\\
\hline
\end{tabular}),
\end{align*}
$
\cL_\cP(\nabla^{\alpha_2}(-2\omega_1+\omega_2))
\simeq
\cL_\cP(\Delta^{\alpha_2}
(-2\omega_1+\omega_2))\in\hat\cE$.
Then, as
\begin{align*}
\hat\cE
&\ni\Delta(\omega_1)
\otimes_\Bbbk\cL_\cP(-3\omega_1)
\simeq
\cL_\cP(\begin{tabular}{|c|}
\hline
$\omega_1$
\\
\hline
$\Delta^{\alpha_2}(-\omega_1+\omega_2)$
\\
\hline
$\varepsilon$
\\
\hline
$\Delta^{\alpha_2}(-2\omega_1+\omega_2)$
\\
\hline
$-\omega_1$
\\
\hline
\end{tabular})
\otimes
(-3\omega_1))
\\
&\simeq
\cL_\cP(\begin{tabular}{|c|}
\hline
$-2\omega_1$
\\
\hline
$\Delta^{\alpha_2}(-4\omega_1+\omega_2)$
\\
\hline
$-3\omega_1$
\\
\hline
$\Delta^{\alpha_2}(-5\omega_1+\omega_2)$
\\
\hline
$-4\omega_1$
\\
\hline
\end{tabular}),
\end{align*}
$
\cL_\cP(\nabla^{\alpha_2}(-5\omega_1+\omega_2))
\simeq
\cL_\cP(\Delta^{\alpha_2}(-5\omega_1+\omega_2))
\in\hat\cE$.
It follows, as
\begin{align*}
\hat\cE
&\ni\Delta(\omega_2)
\otimes_\Bbbk\cL_\cP(-2\omega_1)
\simeq
\cL_\cP(\begin{tabular}{|c|}
\hline
$\Delta^{\alpha_2}(\omega_2)$
\\
\hline
$\omega_1$
\\
\hline
$\Delta^{\alpha_2}(-\omega_1+\omega_2)$
\\
\hline
\begin{tabular}{c|c}
$\Delta^{\alpha_2}(-3\omega_1+2\omega_2)$
&
$\varepsilon$
\end{tabular}
\\
\hline
$\Delta^{\alpha_2}(-2\omega_1+\omega_2)$
\\
\hline
$-\omega_1$
\\
\hline
$\Delta^{\alpha_2}(-3\omega_1+\omega_2)$
\\
\hline
\end{tabular})
\otimes
(-2\omega_1))
\\
&\simeq
\cL_\cP(\begin{tabular}{|c|}
\hline
$\Delta^{\alpha_2}(-2\omega_1+\omega_2)$
\\
\hline
$-\omega_1$
\\
\hline
$\Delta^{\alpha_2}(-3\omega_1+\omega_2)$
\\
\hline
\begin{tabular}{c|c}
$\Delta^{\alpha_2}(-5\omega_1+2\omega_2)$
&
$-2\omega_1$
\end{tabular}
\\
\hline
$\Delta^{\alpha_2}(-4\omega_1+\omega_2)$
\\
\hline
$-3\omega_1$
\\
\hline
$\Delta^{\alpha_2}(-5\omega_1+\omega_2)$
\\
\hline
\end{tabular}),
\end{align*}
that
$
\cL_\cP(\Delta^{\alpha_2}(-5\omega_1+2\omega_2))
\in\hat\cE$.

If $p\geq3$,
$\Delta^{\alpha_2}(-5\omega_1+2\omega_2)
\simeq
\nabla^{\alpha_2}
(-5\omega_1+2\omega_2)$. 
If $p=2$,
as
$\Delta^{\alpha_2}(-5\omega_1+2\omega_2)
=
\begin{tabular}{|c|}
\hline
$L^{\alpha_2}(-5\omega_1+2\omega_2)$
\\
\hline
$-2\omega_1$
\\
\hline
\end{tabular}$
with
$L^{\alpha_2}(-5\omega_1+2\omega_2)$
simple
$P_{\alpha_2}$-module of highest weight 
$-5\omega_1+2\omega_2$,
and as
$\cL_\cP(-2\omega_1)\in\hat\cE$,
so does
$\cL_\cP(\nabla^{\alpha_2}(-5\omega_1+2\omega_2))$.
Thus,
$\cL_\cP(\nabla^{\alpha_2}(-5\omega_1+2\omega_2))\in\hat\cE$
regardless of the characteristic.

Also,
as
\[
\hat\cE\ni\Delta(\omega_1)\otimes_\Bbbk\cO_\cP
\simeq
\cL_\cP(\begin{tabular}{|c|}
\hline
$\omega_1$
\\
\hline
$\Delta^{\alpha_2}(-\omega_1+\omega_2)$
\\
\hline
$\varepsilon$
\\
\hline
$\Delta^{\alpha_2}(-2\omega_1+\omega_2)$
\\
\hline
$-\omega_1$
\\
\hline
\end{tabular}),
\]
one obtains
\[
\cL_\cP(\begin{tabular}{|c|}
\hline
$\omega_1$
\\
\hline
$\nabla^{\alpha_2}(-\omega_1+\omega_2)$
\\
\hline
\end{tabular})
\simeq
\cL_\cP(\begin{tabular}{|c|}
\hline
$\omega_1$
\\
\hline
$\Delta^{\alpha_2}(-\omega_1+\omega_2)$
\\
\hline
\end{tabular})
\in\hat\cE.
\]
Likewise, as
\begin{align*}
\hat\cE
&\ni\Delta(\omega_1)\otimes_\Bbbk\cL(-4\omega_1)
\simeq
\cL_\cP(\begin{tabular}{|c|}
\hline
$\omega_1$
\\
\hline
$\Delta^{\alpha_2}(-\omega_1+\omega_2)$
\\
\hline
$\varepsilon$
\\
\hline
$\Delta^{\alpha_2}(-2\omega_1+\omega_2)$
\\
\hline
$-\omega_1$
\\
\hline
\end{tabular})
\otimes
(-4\omega_1))
\simeq
\cL_\cP(\begin{tabular}{|c|}
\hline
$-3\omega_1$
\\
\hline
$\Delta^{\alpha_2}(-5\omega_1+\omega_2)$
\\
\hline
$-4\omega_1$
\\
\hline
$\Delta^{\alpha_2}(-6\omega_1+\omega_2)$
\\
\hline
$-5\omega_1$
\\
\hline
\end{tabular}),
\end{align*}
one has
\[
\cL_\cP(\begin{tabular}{|c|}
\hline
$\nabla^{\alpha_2}(-6\omega_1+\omega_2)$
\\
\hline
$-5\omega_1$
\\
\hline
\end{tabular})
\simeq
\cL_\cP(\begin{tabular}{|c|}
\hline
$\Delta^{\alpha_2}(-6\omega_1+\omega_2)$
\\
\hline
$-5\omega_1$
\\
\hline
\end{tabular})
\in\hat\cE.
\]
Let
$\tilde\cE$ be the triangulated subcategory of
$\rD^b(\coh\cB)$ Karoubian 
generated by
$\cL_\cB(-r\omega_1)$,
$r\in[0,4]$,
$\cL_\cB(\nabla^{\alpha_2}(-s\omega_1+\omega_2))$
and $\cL_\cB(-s\omega_1+\omega_2)$,
$s\in[2,5]$,
$\cL_\cB(\nabla^{\alpha_2}(-5\omega_1+2\omega_2))$,
$\cL_\cB(-5\omega_1+2\omega_2)$,
$\cL_\cB(\begin{tabular}{|c|}
\hline
$\omega_1$
\\
\hline
$\nabla^{\alpha_2}(-\omega_1+\omega_2)$
\\
\hline
\end{tabular})
$, and
$\cL_\cB(\begin{tabular}{|c|}
\hline
$\nabla^{\alpha_2}(-6\omega_1+\omega_2)$
\\
\hline
$-5\omega_1$
\\
\hline
\end{tabular})
$.
By (1), (2) and
(3) it suffices to show that
all
$\cL_\cB(n\omega_1)\in\tilde\cE$,
$n\in\bbZ$.
As all
$
\cL_\cB(\nabla^{\alpha_2}(-s\omega_1+\omega_2))
\simeq
\cL_\cB(\begin{tabular}{|c|}
\hline
$-s\omega_1+\omega_2$
\\
\hline
$(-s+3)\omega_1-\omega_2$
\\
\hline
\end{tabular})$
and
$\cL_\cB(-s\omega_1+\omega_2)$,
$
s\in[2,5]$,
belong to
$\tilde\cE$,
one obtains
\linebreak
$\cL_\cB((-s+3)\omega_1-\omega_2)\in\tilde\cE$.
Likewise,
as
$\cL_\cB(\nabla^{\alpha_2}(-5\omega_1+2\omega_2))\in\tilde\cE$,
$\cL_\cB(\omega_1-2\omega_2)\in\tilde\cE$.
Then, as
$\nabla(\omega_1)\otimes_\Bbbk\cL_\cB(-\omega_2)\in\tilde\cE$,
$\cL_\cB(2\omega_1-2\omega_2)\in\tilde\cE$.
Likewise, as
$\nabla(\omega_1)\otimes_\Bbbk\cL_\cB(-\rho)\in\tilde\cE$,
$\cL_\cB(-2\omega_2)\in\tilde\cE$.
As
$\nabla(\omega_1)\otimes_\Bbbk
\cL_\cB(-3\omega_1+\omega_2)\in\tilde\cE$,
$\cL_\cB(-4\omega_1+2\omega_2)\in\tilde\cE$.
Then,
as
$\nabla(\omega_2)\otimes_\Bbbk\cL_\cB(-\omega_1)\in\tilde\cE$,
$\cL_\cB(\nabla^{\alpha_2}(-\omega_1+\omega_2))\in\tilde\cE$.
Then,
as
$\cL_\cB(\begin{tabular}{|c|}
\hline
$\omega_1$
\\
\hline
$\nabla^{\alpha_2}(-\omega_1+\omega_2)$
\\
\hline
\end{tabular})
\in\tilde\cE$ by definition,
$\cL_\cB(\omega_1)\in\tilde\cE$.
Likewise,
as
$\nabla(\omega_1)\otimes_\Bbbk
\cL_\cB(-4\omega_1+\omega_2)\in\tilde\cE$,
$\cL_\cB(-6\omega_1+2\omega_2)\in\tilde\cE$.
Then,
as
$\nabla(\omega_2)\otimes_\Bbbk\cL_\cB(-3\omega_1)\in\tilde\cE$,
$\cL_\cB(\nabla^{\alpha_2}(-6\omega_1+\omega_2))\in\tilde\cE$.
Then,
as
$\cL_\cB(\begin{tabular}{|c|}
\hline
$\nabla^{\alpha_2}(-6\omega_1+\omega_2)$
\\
\hline
$-5\omega_1$
\\
\hline
\end{tabular})
\in\tilde\cE$
by definition,
$\cL_\cB(-5\omega_1)\in\tilde\cE$.
Thus
$
\cL_\cB(k\omega_1)\in\tilde\cE$
$\forall k\in[-5,1]$.
As
$\dim\nabla(\omega_1)=7$,
one now obtains all
$\cL_\cB(n\omega_1)\in\tilde\cE$,
$n\in\bbZ$,
from the exact sequence
\begin{multline*}
0\to
\cL_\cB(-7\omega_1)\otimes_\Bbbk
\wedge^7\nabla(\omega_1)
\to
\cL_\cB(-6\omega_1)
\otimes_\Bbbk
\wedge^6\nabla(\omega_1)
\to
\\
\cL_\cB(-5\omega_1)
\otimes_\Bbbk
\wedge^5\nabla(\omega_1)
\to
\dots
\to
\cL_\cB(-\omega_1)
\otimes_\Bbbk
\wedge^1\nabla(\omega_1)
\to
\cO_\cB
\to
0.
\end{multline*}
This finishes a verification in
the case of
$P_{\alpha_2}$
in
type
$\rG_2$.

The other cases are handled entirely similarly.

\setcounter{equation}{0}
\begin{center}
{\bf 
$3^\circ$
Extensions
}
\end{center}

In this section we will compute the extensions among our
$\cE_w$'s given in
(1.7) and (1.9)
to verify (1.10). 

The cases for $G=\SL_2, \SL_3, \Sp_4$
for $P$ a Borel subgroup 
have been done in
\cite{HKR} and \cite{KY}
if $p\geq h$ the Coxeter number of $G$;
the new presentations in (1.7)
somewhat ease the computations in
\cite{KY},
and moreover, adopting those
$\cE_w$'s, $w\in W$, we can get rid of the restrictions on the characteristic.

Let us next explain the characteristic restrictions in type $\rG_2$ stated in 
(1.11).
Let
$P=P_{\alpha_1}$.
One has
\begin{align*}
\Ext^\bullet_\cP
&(\cE_{s_2},
\cE_{s_1s_2})
\simeq
\Ext^\bullet_\cP(\cL_\cP(-\omega_2),
\cL_\cP((-\omega_2)\otimes\begin{tabular}{|c|}
\hline
$
\Delta^{\alpha_1}(2\omega_1-\omega_2)$
\\
\hline
$
\Delta^{\alpha_1}(\omega_1-\omega_2)$
\\
\hline
\end{tabular}
))
\\
&\simeq
\rH^\bullet(\cP, \cL_\cP(
\begin{tabular}{|c|}
\hline
$
\Delta^{\alpha_1}(2\omega_1-\omega_2)$
\\
\hline
$
\nabla^{\alpha_1}(\omega_1-\omega_2)$
\\
\hline
\end{tabular}))
\simeq
\rH^\bullet(\cP, \cL_\cP(
\Delta^{\alpha_1}(2\omega_1-\omega_2)))
\\
&\simeq
\rH^\bullet(\cP, \cL_\cP(
\begin{tabular}{|c|}
\hline
$L^{\alpha_1}(2\omega_1-\omega_2)$
\\
\hline
$\varepsilon$
\\
\hline
\end{tabular}))
\quad\text{as $p=2$}
\end{align*}
with
$L^{\alpha_1}(2\omega_1-\omega_2)$
simple $P_{\alpha_1}$-module of highest weight
$2\omega_1-\omega_2$.
On the other hand, in characteristic $2$
\[
0=\rH^\bullet(\cP, \cL_\cP(
\nabla^{\alpha_1}(2\omega_1-\omega_2)))
=
\rH^\bullet(\cP, \cL_\cP(
\begin{tabular}{|c|}
\hline
$\varepsilon$
\\
\hline
$L^{\alpha_1}(2\omega_1-\omega_2)$
\\
\hline
\end{tabular})),
\]
and hence
$
\rH^i(\cP, \cL_\cP(
L^{\alpha_1}(2\omega_1-\omega_2)))
\simeq
\rH^{i-1}(\cP, \cL_\cP(
\varepsilon))
\simeq
\delta_{i-1,0}\Bbbk
=
\delta_{i1}\Bbbk$,
and (1.11.i) follows.
Likewise (1.11.ii).

Now, to compute all
$\Ext^\bullet_\cB(\cE_x,\cE_y)$,
the case of $G$ in type $\rG_2$ with
$P=B$ is by far the hardest.
As $|W|=12$, there are 144 of them.
Put
$\cB=G/B$.
Let us exhibit the computation
of
$\Ext^\bullet_\cB(\cE_{s_1s_2s_1s_2s_1},\cE_{s_1s_2s_1s_2})$,
which is most complicated among them.
In view of the mal-behaviour in characteristic $2$ and $3$ as noted above,
and in order for
irreducible $G$-module
$\nabla(\omega_1)$ 
not to appear as a composition factor of
$\nabla(\rho)$,
and for other reasons,
we will assume $p\geq7$.
A basic idea is to exploit multiple guises of
the $B$-module structure defining the
$\cE_w$.
For example,
the $B$-module defining
$\cE_{s_1s_2s_1s_2s_1}$
is
$(-\rho)\otimes(\Delta(\omega_2)/(-\omega_2))$
with 
$\Delta(\omega_2)/(-\omega_2)$
having a filtration of
$B$-modules 
\begin{tabular}{|c|}
\hline
$\omega_2$
\\
\hline
$(-\omega_2)\otimes\Delta^{\alpha_1}(3\omega_1)$
\\
\hline
\begin{tabular}{c|c}
$(-\omega_2)\otimes\Delta^{\alpha_1}(2\omega_1)$
&
$\varepsilon$
\end{tabular}
\\
\hline
$(-2\omega_2)\otimes\Delta^{\alpha_1}(3\omega_1)$
\\
\hline
\end{tabular}.
But $\Delta(\omega_2)$
also admits a filtration by
$P_{\alpha_1}$-modules
\[
\Dist(P_{\alpha_1})(\Delta(\omega_2)_{-\alpha_2})
<
\Dist(P_{\alpha_1})(\Bbbk v_2+\Delta(\omega_2)_{-\alpha_2})
<
\Dist(P_{\alpha_1})(\Delta(\omega_2)_{\leq\alpha_1})
\]
with
$v_2$ as in (1.9) and where
$\Delta(\omega_2)_{\leq\alpha_1}
=\sum_{\nu\leq\alpha_1}\Delta(\omega_2)_\nu$
with
$\Dist(P_{\alpha_1})(\Delta(\omega_2)_{-\alpha_2})
=
\begin{tabular}{|c|}
\hline
$\Delta^{\alpha_1}(3\omega_1-2\omega_2)$
\\
\hline
$-\omega_2$
\\
\hline
\end{tabular}$.
Moreover,
$\Delta(\omega_2)/\Dist(P_{\alpha_1})(\Delta(\omega_2)_{-\alpha_2})
=
\begin{tabular}{|c|}
\hline
$\omega_1\otimes\ker(\Delta(\omega_1)\twoheadrightarrow\omega_1)$
\\
\hline
$\Delta^{\alpha_2}(\omega_2)
\otimes(-3\omega_1+\omega_2)$
\\
\hline
$-2\omega_1+\omega_2$
\\
\hline
\end{tabular}$.
Thus
\begin{equation}
\Ext_\cB^\bullet
(\cE_{s_1s_2s_1s_2s_1},\cE_{s_1s_2s_1s_2})
\simeq
\Ext_\cB^\bullet(\begin{tabular}{|c|}
\hline
$\cE_{s_1s_2}$
\\
\hline
$\cL(\Delta^{\alpha_2}(\omega_2)\otimes(-4\omega_1))$
\\
\hline
$\cL(-3\omega_1)$
\\
\hline
$\cL(\Delta^{\alpha_1}(3\omega_1-2\omega_2)\otimes(-\rho))$
\\
\hline
\end{tabular}, \cE_{s_1s_2s_1s_2}).
\end{equation}
By definition
$\cE_{s_1s_2s_1s_2}=\cL_\cB((-\rho)\otimes\{(\Delta(\omega_2)\oplus\Delta(\omega_1))/\Dist(P_{\alpha_2})(\Bbbk(v_4+v_3)+\Delta(\omega_2)_{-3\omega_1+\omega_2})\})$.
We note under the assumption $p\geq7$
that
$\Dist(P_{\alpha_2})(\Bbbk(v_4+v_3)+\Delta(\omega_2)_{-3\omega_1+\omega_2})
=\Dist(P_{\alpha_2})(v_4+v_3)
$,
which admits a $P_{\alpha_2}$-filtration
\begin{tabular}{|c|}
\hline
$(-2\omega_1)\otimes\Delta^{\alpha_2}(\omega_2)$
\\
\hline
$-\omega_1$
\\
\hline
$(-3\omega_1)\otimes
\Delta^{\alpha_2}(\omega_2)$
\\
\hline
\end{tabular}
and also a $B$-filtration
\begin{tabular}{|c|}
\hline
$-2\omega_1+\omega_2$
\\
\hline
$(-\rho)\otimes
\Delta^{\alpha_1}(2\omega_1)$
\\
\hline
$-\omega_2$
\\
\hline
\end{tabular}.

We claim that
$\Ext_\cB^\bullet(\cE_{s_1s_2},
\cE_{s_1s_2s_1s_2})$,
$\Ext_\cB^\bullet(\cL(\Delta^{\alpha_2}(\omega_2)\otimes(-4\omega_1)),\cE_{s_1s_2s_1s_2})$
\linebreak
and 
$\Ext_\cB^\bullet(\cL(-3\omega_1),
\cE_{s_1s_2s_1s_2}))$ all vanish.
First,
$
\Ext_\cB^\bullet(\cE_{s_1s_2}, \cE_{s_1s_2s_1s_2})
\simeq
\Ext_\cB^\bullet(\cL((-\omega_2)\otimes
\ker(\Delta(\omega_1)\twoheadrightarrow\omega_1)),\cE_{s_1s_2s_1s_2})
$
with
$\forall i\in\bbN$,
\begin{align*}
\Ext_\cB^i
&(\cL(-\omega_2+\omega_1),
\cE_{s_1s_2s_1s_2})
\\
&\simeq
\rH^i(\cB, 
\cL((\omega_2-\omega_1-\rho)\otimes\{
(\Delta(\omega_2)\oplus\Delta(\omega_1))/\Dist(P_{\alpha_2})(v_4+v_3)\}))
\\
&\simeq
\rH^i(\cB, 
\cL((-2\omega_1)\otimes\{
(\Delta(\omega_2)\oplus\Delta(\omega_1))/\Dist(P_{\alpha_2})(v_4+v_3)\}))
\\
&\simeq
\rH^{i+1}(\cB, 
\cL((-2\omega_1)\otimes\Dist(P_{\alpha_2})(v_4+v_3)))
\\
&\simeq
\rH^{i+1}(\cB, 
\cL((-2\omega_1)\otimes\begin{tabular}{|c|}
\hline
$(-2\omega_1)\otimes\Delta^{\alpha_2}(\omega_2)$
\\
\hline
$-\omega_1$
\\
\hline
$(-3\omega_1)\otimes
\Delta^{\alpha_2}(\omega_2)$
\\
\hline
\end{tabular}
))
\simeq
\rH^{i+1}(\cB, 
\cL(\begin{tabular}{|c|}
\hline
$\Delta^{\alpha_2}(-4\omega_1+\omega_2)$
\\
\hline
$-3\omega_1$
\\
\hline
$
\Delta^{\alpha_2}(-5\omega_1+\omega_2)$
\\
\hline
\end{tabular}
))
\\
&=0
\quad\text{as $p\geq5$}
\end{align*}
while
$
\Ext_\cB^i(\cL((-\omega_2)\otimes
\Delta(\omega_1)),\cE_{s_1s_2s_1s_2})
\simeq
\Ext_\cB^i(\cE_{s_2}, \cE_{s_1s_2s_1s_2})
\otimes\Delta(\omega_1)^*
$
with
\begin{align*}
\Ext_\cB^i
&(\cE_{s_2}, \cE_{s_1s_2s_1s_2})
\\
&\simeq
\Ext_\cB^i(\cL(-\omega_2), \cL((-\rho)\otimes\{
(\Delta(\omega_2)\oplus\Delta(\omega_1))/\Dist(P_{\alpha_2})(v_4+v_3)\}))
\\
&\simeq
\rH^i(\cB, 
\cL((-\omega_1)\otimes\{
(\Delta(\omega_2)\oplus\Delta(\omega_1))/\Dist(P_{\alpha_2})(v_4+v_3)\}))
\\
&\simeq
\rH^{i+1}(\cB, 
\cL((-\omega_1)\otimes\Dist(P_{\alpha_2})(v_4+v_3)))
\\
&\simeq
\rH^{i+1}(\cB, 
\cL((-\omega_1)\otimes\begin{tabular}{|c|}
\hline
$(-2\omega_1)\otimes\Delta^{\alpha_2}(\omega_2)$
\\
\hline
$-\omega_1$
\\
\hline
$(-3\omega_1)\otimes
\Delta^{\alpha_2}(\omega_2)$
\\
\hline
\end{tabular}
))
\simeq
\rH^{i+1}(\cB, 
\cL(\begin{tabular}{|c|}
\hline
$\Delta^{\alpha_2}(-3\omega_1+\omega_2)$
\\
\hline
$-2\omega_1$
\\
\hline
$
\Delta^{\alpha_2}(-4\omega_1+\omega_2)$
\\
\hline
\end{tabular}
))
\\
&=0
\quad\text{by \cite[II.6.18]{J}}.
\end{align*}
and hence
$\Ext_\cB^\bullet(\cE_{s_1s_2}, \cE_{s_1s_2s_1s_2})
=0$.
Next,
\begin{align*}
\Ext^\bullet_\cB
&(\cL(\Delta^{\alpha_2}(\omega_2)
\otimes(-4\omega_1)),
\cE_{s_1s_2s_1s_2})
\simeq
\Ext^\bullet_\cB
(\cL(\Delta^{\alpha_2}(\omega_2)
\otimes(-4\omega_1)),
\\
&\hspace{2cm}
\cL((-\rho)\otimes\{
(\Delta(\omega_2)\oplus\Delta(\omega_1))/\Dist(P_{\alpha_2})(v_4+v_3)\}))
\\
&\simeq
\rH^\bullet(\cB,
\cL((-\omega_2)\otimes3\omega_1\otimes
\nabla^{\alpha_2}(-3\omega_1+\omega_2)
\otimes
\\
&\hspace{2cm}
\{
(\Delta(\omega_2)\oplus\Delta(\omega_1))/\Dist(P_{\alpha_2})(v_4+v_3)\}))
=0
\end{align*}
as
$3\omega_1\otimes
\nabla^{\alpha_2}(-3\omega_1+\omega_2)
\otimes
\{
(\Delta(\omega_2)\oplus\Delta(\omega_1))/\Dist(P_{\alpha_2})(v_4+v_3)\}$
is equipped with a structure of
$P_{\alpha_2}$-module.
Likewise,
\begin{align*}
\Ext_\cB^i
&(\cL(-3\omega_1), \cE_{s_1s_2s_1s_2})
\\
&\simeq
\Ext_\cB^i(\cL(-3\omega_1), 
\cL((-\rho)\otimes\{
(\Delta(\omega_2)\oplus\Delta(\omega_1))/\Dist(P_{\alpha_2})
(v_4+v_3)\}))
\\
&\simeq
\rH^i(\cB,
\cL((-\omega_2)\otimes2\omega_1\otimes\{
(\Delta(\omega_2)\oplus\Delta(\omega_1))/\Dist(P_{\alpha_2})(v_4+v_3)\}))
=0.
\end{align*}

It now follows 
$\forall i\in\bbN$ that
\begin{align*}
\Ext_\cB^i
&(\cE_{s_1s_2s_1s_2s_1}, 
\cE_{s_1s_2s_1s_2})
\simeq
\Ext_\cB^i(\cL(\Delta^{\alpha_1}(3\omega_1-2\omega_2)\otimes(-\rho)),
\cE_{s_1s_2s_1s_2})
\\
&\simeq
\Ext_\cB^i(\cL(\Delta^{\alpha_1}(3\omega_1-2\omega_2)\otimes(-\rho)),
\\
&\hspace{2cm}
\cL((-\rho)\otimes\{
(\Delta(\omega_2)\oplus\Delta(\omega_1))/\Dist(P_{\alpha_2})(v_4+v_3)\}))
\\
&\simeq
\rH^i(\cB,
\cL(\nabla^{\alpha_1}(3\omega_1-\omega_2)\otimes
\{
(\Delta(\omega_2)\oplus\Delta(\omega_1))/\Dist(P_{\alpha_2})(v_4+v_3)\}))
\\
&\simeq
\rH^{i+1}(\cB,
\cL(\nabla^{\alpha_1}(3\omega_1-\omega_2)\otimes
\Dist(P_{\alpha_2})(v_4+v_3)))
\\
&\simeq
\rH^{i+1}(\cB,
\cL(\nabla^{\alpha_1}(3\omega_1-\omega_2)\otimes
\begin{tabular}{|c|}
\hline
$-2\omega_1+\omega_2$
\\
\hline
$(-\rho)\otimes
\Delta^{\alpha_1}(2\omega_1)$
\\
\hline
$-\omega_2$
\\
\hline
\end{tabular}
))
\\
&\simeq
\rH^{i+1}(\cB,
\cL(
\begin{tabular}{|c|}
\hline
$\nabla^{\alpha_1}(3\omega_1)\otimes(-2\omega_1)$
\\
\hline
$\nabla^{\alpha_1}(3\omega_1-\omega_2)\otimes(-\omega_1)\otimes
\nabla^{\alpha_1}(2\omega_1-\omega_2)$
\\
\hline
$\nabla^{\alpha_1}(3\omega_1-2\omega_2)$
\\
\hline
\end{tabular}
))
\end{align*}
with
\begin{align*}
\rH^{i+1}(\cB,
&\cL(\nabla^{\alpha_1}(3\omega_1)\otimes(-2\omega_1)
))
\simeq
\rH^{i}(\cB,
\cL(\nabla^{\alpha_1}(3\omega_1)\otimes(-\omega_2)
))
\quad\text{by the Serre duality
}
\\
&\simeq
\rH^{i}(\cB,
\cL(\nabla^{\alpha_1}(3\omega_1-\omega_2)
))
\\
&
=0
=
\rH^\bullet(\cB,
\cL(\nabla^{\alpha_1}(3\omega_1-\omega_2)\otimes(-\omega_1)\otimes
\nabla^{\alpha_1}(2\omega_1-\omega_2)
)).
\end{align*}
One thus obtains that
\begin{align*}
\Ext_\cB^i(\cE_{s_1s_2s_1s_2s_1}, \cE_{s_1s_2s_1s_2})
&\simeq
\rH^{i+1}(\cB,
\cL(\nabla^{\alpha_1}(3\omega_1-2\omega_2)
))
\\
&\simeq
\delta_{i+1,1}\Bbbk
=\delta_{i0}\Bbbk
\quad\text{by \cite[II.6.18]{J}}.
\end{align*}

\setcounter{equation}{0}
\begin{center}
$4^\circ$
{\bf
Kapranov's sheaves
}
\end{center}

In \cite{K08}/\cite{KNS} we showed that Kapranov's 
sheaves on the Grassmannians 
\cite{Kap}
constitute a  tilting sheaf in positive characteristic if 
the characteristic is large enough,
and their parametrization by
$W^P$ verifies Catanese's conjecture.
In this section 
We will briefly discuss 
Kapranov's sheaves on the flag variety of
$\GL_3$
and on the quadrics
for future study.

\setcounter{equation}{0}
\noindent
(5.1)
Let us first consider the flag variety
$\cB=\GL(E)/B$
with
$E$ of dimension 3.
If $p\geq3$, as well as in characteristic 0,
Kapranov's sheaves
$\cE_e
=
\cL_\cB(\nabla^{{\alpha_2}}(-\omega_1)\otimes(-\omega_2))
\simeq
\cL_\cB(-\rho)$,
$\cE_{s_1}
=
\cL_\cB(\nabla^{{\alpha_2}}(-\omega_1+\omega_2)\otimes(-\omega_2))
\simeq
\cL_\cB(\nabla^{{\alpha_2}}(\omega_2)
\otimes(-\rho))$,
$\cE_{s_2s_1}
=
\cL_\cB(-\omega_2)$,
$\cE_{s_2}
=
\cL_\cB(\nabla^{{\alpha_2}}(-\omega_1))
\simeq
\cL_\cB(-\omega_1)$,
$\cE_{s_1s_2}
=
\cL_\cB(\nabla^{{\alpha_2}}(-\omega_1+\omega_2))
\simeq
\cL_\cB(\nabla^{{\alpha_2}}(\omega_2)
\otimes(-\omega_1))$,
and
$\cE_{s_2s_1s_2}
=
\cO_\cB
$
from
\cite{Kap88}/\cite[2.2.3]{Boh}
form a complete strongly 
exceptional sequence
on the flag variety.
The nonvanishing
$\Mod_\cB(\cE_x, \cE_y)$,
$x\ne y$, are,
however, given by
\begin{alignat*}{2}
\Mod_\cB(\cE_e, \cE_{s_1})
&\simeq
\nabla(\omega_2),
\quad&
\Mod_\cB(\cE_e, \cE_{s_2})
&\simeq
\nabla(\omega_2),
\\
\Mod_\cB(\cE_e, \cE_{s_1s_2})
&\simeq
\begin{tabular}{|c|}
\hline
$\nabla(2\omega_2)$
\\
\hline
$\nabla(\omega_1)$
\\
\hline
\end{tabular},
\quad&
\Mod_\cB(\cE_e, \cE_{s_2s_1})
&\simeq
\nabla(\omega_1),
\\
\Mod_\cB(\cE_e, \cE_{s_2s_1s_2})
&\simeq
\nabla(\rho),
\quad&
\Mod_\cB(\cE_{s_1}, \cE_{s_2})
&\simeq
\Bbbk,
\\
\Mod_\cB(\cE_{s_1}, \cE_{s_2s_1})
&\simeq
\nabla(\omega_2),
\quad&
\Mod_\cB(\cE_{s_1}, \cE_{s_1s_2})
&\simeq
\nabla(\omega_2)^{\oplus2},
\\
\Mod_\cB(\cE_{s_1}, \cE_{s_2s_1s_2})
&\simeq
\begin{tabular}{|c|}
\hline
$\nabla(2\omega_2)$
\\
\hline
$\nabla(\omega_1)$
\\
\hline
\end{tabular},
\qquad&
\Mod_\cB(\cE_{s_2}, \cE_{s_1s_2})
&\simeq
\nabla(\omega_2),
\end{alignat*}\begin{alignat*}{2}
\Mod_\cB(\cE_{s_2}, \cE_{s_2s_1s_2})
&\simeq
\nabla(\omega_1),
\quad&
\Mod_\cB(\cE_{s_1s_2}, \cE_{s_2s_1s_2})
&\simeq
\nabla(\omega_2),
\\
\Mod_\cB(\cE_{s_2s_1}, \cE_{s_1s_2})
&\simeq
\Bbbk,
\quad&
\Mod_\cB(\cE_{s_2s_1}, \cE_{s_2s_1s_2})
&\simeq
\nabla(\omega_2).
\end{alignat*}
In particular,
$|\{w\in W\setminus
\{s_1s_2\}\mid
\Mod_\cB(\cE_{w}, \cE_{s_1s_2})\ne0\}|=4$,
and hence
there is no reindexing of these $\cE_w$
by $W$ such that the Catanese conjecture
hold,
contrary to our construction
(1.8), (1.10).

\setcounter{equation}{0}
\noindent
(5.2)
Let us next consider the quadric
$\cQ=\cQ_n$ of 
dimension $n\geq3$ in odd characteristic.
Let
$E$ be a $\Bbbk$-linear space of dimension $n+2$.
In case $n$ is odd, write $n=2m+1$.
Let
$e_1, \dots, e_{m+1}, e_0, e_{-m-1},
\dots,e_{-1}$ 
be a basis of $E$ and define a quadratic form $q$ on $E$ by
$q(\sum_{k=-m-1}^{m+1}x_ke_k)=
x_1x_{-1}+\dots+x_{m+1}x_{-m-1}+x_0^2$.
If
$n$ is even, write $n=2m$.
Let
$e_1, \dots, e_{m+1}, e_{-m-1},
\dots,e_{-1}$ 
be a basis of $E$ equipped with quadratic form
$q(\sum_{k=1}^{m+1}(x_ke_k+x_{-k}e_{-k})=
x_1x_{-1}+\dots+x_{m+1}x_{-m-1}$.
In either case let
$G=\mathrm{SO}(E)$,
$P=\rN_G(\Bbbk e_{-1})$.
Then a closed immersion
$i : G/P\to\bbP(E)$
via $gP\mapsto[ge_{-1}]$
identifies $G/P$
with
$\cQ$.

If $n=2m+1$, 
$T=\{\diag(\zeta_1,\dots,\zeta_{m+1},1,\zeta_{m+1}^{-1},\dots,\zeta_{1}^{-1})\mid\zeta_i\in\Bbbk^\times
\ \forall i\}$
forms a maximal torus of $G$.
Take as simple roots
$
\alpha_1=
\varepsilon_1-\varepsilon_2,
\dots,
\alpha_m=
\varepsilon_m-\varepsilon_{m+1},
\alpha_{m+1}=\varepsilon_{m+1}$
with
$\varepsilon_k : \diag(\zeta_1,
\dots,\zeta_{m+1},\zeta_0,\zeta_{-m-1},\dots,\zeta_{-1})
\mapsto\zeta_k$.
Then the fundamental weights are given by
$
\omega_1=\varepsilon_1,
\omega_2=\varepsilon_1+\varepsilon_2,
\dots,
\omega_m=
\varepsilon_1+\dots+\varepsilon_m,
\omega_{m+1}=
\frac{1}{2}(\varepsilon_1+\dots+\varepsilon_{m+1})$,
and
$W^P=\{e, s_1, s_2s_1, \dots, s_{m+1}s_m\dots
s_2s_1, s_ms_{m+1}s_m\dots s_2s_1,
s_{m-1}s_ms_{m+1}s_m\dots s_2s_1,
\dots,
\linebreak
s_2\dots s_{m-1}s_ms_{m+1}s_m\dots s_2s_1,
s_1s_2\dots s_{m-1}s_ms_{m+1}s_m\dots s_2s_1\}$.
Define
$\cE_e=\cO_\cQ$,
$\cE_{s_1}=
\linebreak
\cO_\cQ(-1)$,
$\cE_{s_2s_1}=\cO_\cQ(-2)$,
\dots,
$\cE_{s_m\dots s_2s_1}=\cO_\cQ(-m)$,
$\cE_{s_{m+1}s_m\dots s_2s_1}
=
\linebreak
\cL(\nabla^P(\omega_{m+1}))(-m-1)$,
$\cE_{s_ms_{m+1}s_m\dots s_2s_1}=\cO_\cQ(-m-1)$,
\dots,
$\cE_{s_2\dots
s_ms_{m+1}s_m\dots s_2s_1}=
\linebreak
\cO_\cQ(-n+2)$,
and 
$\cE_{s_1\dots
s_ms_{m+1}s_m\dots s_2s_1}=\cO_\cQ(-n+1)$.

If $n=2m$, 
$T=\{\diag(\zeta_1,\dots,\zeta_{m+1},\zeta_{m+1}^{-1},\dots,\zeta_{1}^{-1})\mid\zeta_i\in\Bbbk^\times
\ \forall i\}$
forms a maximal torus of $G$.
Take as simple roots
$
\alpha_1=
\varepsilon_1-\varepsilon_2,
\dots,
\alpha_{m-1}=
\varepsilon_{m-1}-\varepsilon_{m},
\alpha_m=
\varepsilon_m-\varepsilon_{m+1},
\alpha_{m+1}=\varepsilon_m+\varepsilon_{m+1}$
with
$\varepsilon_k : \diag(\zeta_1,
\dots,\zeta_{m+1},\zeta_{-m-1},\dots,\zeta_{-1})
\mapsto\zeta_k$.
Then the fundamental weights are $
\omega_1=\varepsilon_1,
\omega_2=\varepsilon_1+\varepsilon_2,
\dots,
\omega_{m-1}=
\varepsilon_1+\dots+\varepsilon_{m-1},\omega_m=
\frac{1}{2}(
\varepsilon_1+\dots+\varepsilon_m-\varepsilon_{m+1}),
\omega_{m+1}=
\frac{1}{2}(\varepsilon_1+\dots+\varepsilon_m+\varepsilon_{m+1})$,
and
$W^P=\{e, s_1, s_2s_1, \dots, s_ms_{m-1}\dots
s_1, s_{m+1}s_{m-1}\dots s_2s_1,
\linebreak
s_ms_{m+1}s_{m-1}\dots s_1,
s_{m-1}s_ms_{m+1}s_{m-1}\dots s_1,
s_{m-2}s_{m-1}s_ms_{m+1}s_{m-1}\dots s_1,
\dots,
\linebreak
s_1s_2\dots s_{m-1}s_ms_{m+1}s_{m-1}\dots s_1\}$.
Define
$\cE_e=\cO_\cQ$,
$\cE_{s_1}=\cO_\cQ(-1)$,
$\cE_{s_2s_1}=\cO_\cQ(-2)$,
\dots,
$\cE_{s_{m-1}\dots s_2s_1}=\cO_\cQ(-(m-1))$,
$\cE_{s_m\dots s_2s_1}=\cL_\cQ(\nabla^P(\omega_{m+1}))(-m)$,
$\cE_{s_{m+1}s_{m-1}\dots s_2s_1}=
\linebreak
\cL(\nabla^P(\omega_{m}))(-m)$,
$\cE_{s_ms_{m+1}s_{m-1}\dots s_2s_1}
=
\cO_\cQ(-m)$,
\dots,
$\cE_{s_2\dots
s_ms_{m+1}s_{m-1}\dots s_2s_1}=\cO_\cQ(-n+2)$,
and 
$\cE_{s_1\dots
s_ms_{m+1}s_m\dots s_2s_1}=\cO_\cQ(-n+1)$.

In either case
Langer
\cite{La}
shows that
the 
$\cE_w$,
$w\in W^P$,
verify Catanese's conjecture,
and that
for
$p\geq n+1$ the Coxeter number of
$G$
all $\cE_w$'s appear as direct summands
of
$F_*\cO_\cQ$.
More precisely,
let
$A=\rS_\Bbbk(E^*)/(q)$
and put
$\bar A=A/(a^p\mid a\in E^*)$.
Let
$\bar A_j$ be the $j$-th homogeneous part of
$\bar A$.
If $n=2m+1$,
\begin{multline*}
F_*\cO_\cQ
\simeq
\{
\coprod_{i\in[0,n[\setminus\{m+1\}}\cO_Q(-i)\otimes_\Bbbk\bar A_{ip}
\}
\oplus
\{
\cO_\cQ(-m-1)\otimes_\Bbbk
\bar A_{mp-n}
\}
\oplus
\\
\cL(\nabla^P(\omega_{m+1}))^{\oplus  r}
(-m-1)
\end{multline*}
with
$r=\frac{
\dim\bar A_{(m+1)p}-
\dim\bar A_{mp-n}}{\dim L(\omega_{m+1})
}$,
while for $n=2m$
\begin{multline*}
F_*\cO_\cQ
\simeq
\{
\coprod_{i\in[0,n[\setminus\{m\}}\cO_Q(-i)\otimes_\Bbbk\bar A_{ip}
\}
\oplus
\{
\cO_\cQ(-m)\otimes_\Bbbk
\bar A_{mp-n}
\}
\oplus
\\
\{\cL(\nabla^P(\omega_{m})\oplus
\nabla^P(\omega_{m+1}))\}^{\oplus  s}(-m)
\end{multline*}
with
$s=\frac{
\dim\bar A_{mp}-
\dim\bar A_{mp-n}}{\dim \{
L(\omega_{m})\oplus
L(\omega_{m+1})\}
}
$.
Note that our parametrization of the
$\cE_w$ is different from
that of
B\"{o}hning's, but is consistent with
(1.10) in case $n=3$.

\end{document}